\documentclass[12pt]{article}
\usepackage{graphicx}
\usepackage{picture}
\usepackage{hyperref}
\usepackage{color}
\usepackage[english]{babel}
\usepackage{amsmath,amsthm}
\usepackage{amssymb}
\usepackage{enumerate}
\usepackage{amsfonts}
\usepackage{a4wide}
\usepackage{layout}
\usepackage{latexsym}
\usepackage{amsthm}
\numberwithin{equation}{section}
\theoremstyle{plain}
\newtheorem{theorem}{Theorem}[section]
\newtheorem{lemma}[theorem]{Lemma}
\newtheorem{corollary}[theorem]{Corollary}

\newtheorem{proposition}[theorem]{Proposition}
\newcommand{\eps}{\varepsilon}

\newcommand{\nn}{\nonumber}
\newcommand{\dis}{\displaystyle}
\newcommand{\und}{\underline}

\newcommand{\la}{\lambda}

\newcommand{\ga}{\gamma}
\newcommand{\Ga}{\Gamma}

\newcommand{\si}{\sigma}

\begin{document}

\title{Non local branching Brownians with annihilation\\
 and free boundary problems}

\author{Anna De Masi \and Pablo A. Ferrari \and Errico Presutti \and Nahuel Soprano-Loto}

%
%

\date{}

\maketitle
\parindent 0pt

\abstract{We study a system of branching Brownian motions on $\mathbb R$ with
 annihilation: at each branching time  a new particle is created and the leftmost one  is deleted. In \cite{DFPS} it has been studied the case of strictly local creations (the new particle is put exactly at the same position of the branching particle), in \cite{DR11} instead the position $y$ of the new particle   has a distribution $p(x,y)dy$, $x$ the position of the branching particle, however particles in between branching times do not move.  In this paper we consider Brownian motions as in \cite{DFPS} and non local branching as in  \cite{DR11} and prove convergence in the continuum limit (when the number $N$ of particles diverges) to  a limit density which satisfies a free boundary problem when this has classical solutions, local in time existence of classical solution has been proved recently in \cite{JM2}.  We use in the convergence a stronger topology than in \cite{DFPS} and  \cite{DR11} and have explicit bounds on the rate of convergence.}

\paragraph{Keywords} Hydrodynamic limit. Free boundary problems. Branching Brownian motion. Brunet-Derrida models.


\section{Introduction}

The system considered in this  paper fits in a class of models proposed by
Brunet and Derrida  in \cite{BD97} to study selection mechanisms in biological systems and continues a line of research initiated by
Durrett and Remenik in \cite{DR11}.

Durrett and Remenik
have in fact studied a model of  particles on $\mathbb R$ which
independently at rate 1 branch creating a new particle whose position is chosen randomly with probability $p(x,y)dy$, $p(x,y)=p(0,y-x)$, if $x$ is the position of the generating particle.
At the same time the leftmost particle is deleted so that the total number of particles is constant.

In the
biological interpretation  particles are individuals in a population, the position of a particle is  ``its degree of fitness'', the larger the position the higher the fitness. If the environment  supports only populations of a given size then to each birth there must correspond a death. The removal of the leftmost and hence less fitted  particle is a very effective Darwinian selection rule to implement the conservation of  the  population size. Even if the  duplication rule is regressive, i.e.\ $p(x,y)$ has support on $y<x$, nonetheless the population fitness improves and if $p(0,x)>0$ for $x\in (-a,0)$, $a>0$,    then as    time diverges the whole population concentrates around  the position of the initially best fitted individual.  Durrett and  Remenik have studied the case where
$p(x,y)$ is symmetric and discussed the occurrence of traveling solutions which describe a steady improvement of the population fitness, see also \cite{Mai16} and Brunet, Derrida \cite{BD97}, \cite{BBD17} for the analysis of traveling waves in a large class of systems.

Such issues are better studied in the continuum  limit $N\to \infty$, the natural guess for the continuum version of the Durrett and  Remenik duplication process is
 \begin{equation}
 \label{1.1}
 \frac{\partial}{\partial t} \rho(x,t)= \int_{X_t}^\infty dy \; p(y,x) \rho(y,t)dy,\; x \ge X_t; \quad \rho(r,0)=\rho_0(r)
 \end{equation}
where $\rho(x,t)$ is the particles density and $X_t = \sup \{r:\rho(r,t) =0\}$.  The removal process in the continuum is more implicit and given by the condition on $X_t$ that
for all $t\ge 0$
 \begin{equation}
 \label{1.2}
 \int_{X_t}^\infty dy \;   \rho(y,t) =1
 \end{equation}
Under suitable assumptions on the initial datum $\rho_0$ and on the probability kernel $p(x,y)$
Durrett and Remenik have  proved that  \eqref{1.1} and
\eqref{1.2} have a unique solution and that this is the limit density of the particles system.

\eqref{1.1} has a nice probabilistic interpretation:
\begin{equation}
 \label{1.3}
  \rho(x,t)= e^t E_x\Big[\rho_0(x^*(t)) \mathbf 1_{\tau >t}\Big]
 \end{equation}
where $x^*(t)$ is the jump Markov process with generator
\begin{equation}
 \label{1.4}
A^* f(x)  = \int dy\, p^*(x,y) \Big( f(y)-f(x)\Big),\quad  p^*(x,y)= p(y,x)
 \end{equation}
and
\begin{equation}
 \label{1.5}
 \tau = \inf\{t: x^*(t) \le X_t\}
 \end{equation}
\eqref{1.3} is the backward Kolmogorov equation for the jump process with generator $A$ equal to the adjoint of $A^*$ (i.e.\ when the jump $x \to y$ has probability $p(x,y)$).

Similar formulas hold as well in the case considered in this paper  where we study
a natural extension of the Durrett-Remenik model   where   particles move  as independent  Brownian motions   in between branching times: biologically this means that the individuals fitness changes randomly in time. As in \cite{DR11}
the new particles are created at random positions with  probability $p(x,y)$ and, like before,  as soon as a particle is created the leftmost one is deleted.  We have already studied in a previous paper, \cite{DFPS}, the case where $p(x,y)=\delta(y-x)$, namely when the duplication is exact.  The extension to  general $p(x,y)$, which is the aim of this paper, is obtained following the same scheme used in  \cite{DFPS}, and indeed some of the proofs are straightforward adaptation of those in  \cite{DFPS} and its details are omitted.   We thus focus on the new parts for which we give complete details.

The main novelty in our case is that the continuum version of the model gives rise to a new free boundary problem  which does not seem to fit in the class already studied in the literature.  The main difficulty is that the motion of the edge is not monotonous as in
 \cite{DR11} and cannot be reduced to a classical Stefan problem as in  \cite{DFPS} by studying the spatial derivative of the limit density.  The other new aspect of the present paper is the use of a stronger topology in the convergence of the process to a limit density which is related to mass transport in the comparison of probability measures.

In the next two sections we make precise the model  and state the main results, an outline of how the paper is organized is then given at the end of  Section \ref{sec:DD2.1}.
We conclude this introduction by mentioning that there have been several papers about particles processes which in the continuum limit are described by
free boundary problems.  Some of them will be mentioned in the sequel, for a list we refer to a survey   on the subject, \cite{CDGP}.

\section{The model}
\label{sec:DD2}

We will consider in this paper several processes, the main one is $\und x(t)$ (also called sometimes the ``true process'').  We will also use auxiliary processes $\und x^{\delta,\pm}(t)$, $\delta$ a positive parameter,  called the upper and lower stochastic barriers.  We define all of them in a same space as subsets of a  ``basic'' process $\und y(t)$ that we define next.

\smallskip

{\it {The ``basic'' process $\und y(t)$}}.

\noindent The state space of the ``basic'' process $\und y(t)$ is the set of configurations with finitely many point particles, we will denote by $|\und y|$ the number of particles in  a configuration $\und y$. If convenient we may label the particles by writing
$\und y(t)=(y_1(t),\dots,y_n(t))$, sequences which only differ for the labelling are however considered equivalent.

To define the process we attach to each particle  an independent exponential clock of intensity 1: when it rings for a particle (call $x$ its   position when its clock rings) then a new particle is created at position $y$ with probability $p(x,y)dy$, assumptions on $p(x,y)$ are stated later.  In between clock rings   the  particles  are independent
Brownian motions.  $P$ denotes the law of the process and when needed we will write
$P^{\und y_0}$ to specify that the process starts from $\und y_0$.

\smallskip

{\it {The ``true'' process $\und x(t)$}}.

\noindent
As in the basic process each particle of $\und x(t)$ has an independent exponential clock of intensity 1 and in between clock rings   the  particles  are independent
Brownian motions.  Also for $\und x(t)$  when a clock  rings for a particle (call $x$ its   position) then (like in the basic process) a new particle is created at position $y$ with probability $p(x,y)dy$, here however
at the same time when a new particle is created   the leftmost particle (among those previously present plus the new one) is deleted, so that $|\und x(t)|$ is constant.
Evidently  $\und x(t)$ can be realized as a subset of  $\und y(t)$ obtained by disregarding in $\und y(t)$ the  particles which in  $\und x(t)$ are deleted as well as   all their descendants.

\smallskip

{\it {The stochastic barriers  $\und x^{\delta,\pm}(t)$}}.

\nopagebreak
\noindent They are defined  like $\und x(t)$, (see Subsection \ref{subsec2.1}) the difference being that the removal of particles is not simultaneous to the branching but it occurs at discrete times $k\delta$, $k\in \mathbb N$, $\delta>0$.

\smallskip

{\it {The initial configuration}}.

\noindent We will study the process $\und x(t)$ having fixed (arbitrarily) the initial number of particles, denoted by  $N$.  The process starts from an initial configuration $\und x_0$ whose distribution is obtained by taking $N$ independent copies of a position variable distributed with probability $\rho_0(r) dr$. We suppose that $\rho_0(r)$ is a smooth ($C^\infty$) probability density with support in $[-A,A]$, $A>0$. Thus $\und x(t)$ is realized as a subset of the basic process which starts from $\und y(0)=\und x_0$ and its law will be denoted by $P^{(N)}$.

\smallskip

{\it {Assumptions on  $p(x,y)$}}.

\noindent The jump probability density  $p(x,y)$ is translation invariant, namely $p(x,y)=p(0,y-x)$. $p(0,x)$ is a smooth ($C^\infty$) function with finite range:
$p(0,x) =0$ when $|x| \ge \xi$, $\xi>0$.

\smallskip

The assumptions on $\rho_0$ and $p(x,y)$ could be relaxed but we would have some more technical details to take care of.

\smallskip

{\it The counting measure}.

\noindent Given a particle configuration $\und x$ we call
	 \begin{equation}
	 \label{DD2.1}
\pi_{\und x}(dr) := \sum_{x\in  \und x} \delta_{x}(r)dr
   \end{equation}
 the associated counting measure.  Our aim is to study the behavior of the probability measure
$\frac 1N\pi_{\und x(t)}(dr)$, $N= |\und x|$, for large $N$ in a finite time interval $[0,T]$, $T>0$.

\smallskip

{\it {Mass transport order}}.

\noindent Let $\mu$ and $\nu$ be finite, positive measures on $\mathbb R$, we then write for $\eps\ge 0$;
	\begin{align}
	\label{DD2.2}
\mu \preccurlyeq \nu \;\, \text{modulo $\eps$ if for all $r$}
\;\; \int_r^\infty \mu(dr') \le \int_r^\infty  \nu(dr')+ \eps
    \end{align}
When $\eps=0$ we just write $\mu \preccurlyeq \nu$.

\smallskip

If $\mu \preccurlyeq \nu$ we can obtain $\nu$ from $\mu$ by moving some of the mass of $\mu$ to the right, hence the name ``mass transport order''.  This notion will be used in the next section to define a topology (in the mass transport order) and in fact our main result, stated in Theorem \ref{thmDD2.1}, is that the counting measure of the true process is close in the sense of mass transport order to a probability density $u(r,t)$.

\section{Main results}
\label{sec:DD2.1}
As mentioned in the previous section our main result is that the counting measure of the true process is close in the sense of mass transport order to a probability density $u(r,t)$.   To quantify the notion we first define a distance between two probability measures $\mu$ and $\nu$ on $\mathbb R$ as
	\begin{equation}
	\label{DD2.4.0}
|\mu-\nu|_{\rm mt}:= \sup_{r\in \mathbb R}\Big| \mu\Big[[r,\infty)\Big]-\nu\Big[[r,\infty)\Big]\Big|
\Big\}
     \end{equation}
We then set:

\smallskip

{\it Neighborhoods in the mass transport order}.

\noindent
Let $\mathcal M$ be the space of probability measures valued functions on $\mathbb R_+$ whose elements are denoted by $(\mu_t)_{t\ge 0}$.
Let $u(r,t)$ be for each $t\ge 0$ a probability density on $\mathbb R$ so that $u(r,t)dr \in \mathcal M$.
We define below neighborhoods $\mathcal X_{T,\eps,n}(u)$ in $\mathcal M$ of $u$ (in the sense of
mass transport order), where $T>0$, $\eps>0$, $n \in \mathbb N$.
 Calling
     $$
     t_k= k 2^{-n}T,\quad k=0,1,..,2^n
     $$
     we set
     	\begin{equation}
	\label{DD2.4.1}
\mathcal X_{T,\eps,n}(u) = \bigcap_{k=0}^{2^n} \Big\{
|\mu_{t_k}(dr)- u(r,t_k)dr|_{\rm mt} <\eps\Big\}
     \end{equation}
 This means
 	\begin{align}
	\label{DD2.4bis}
\mathcal X_{T,\eps,n}(u) =  \Big\{\sup_{r\in \mathbb R}\;\;\max_{ k=0,1,..,2^n}\Big |\int_r^\infty
\mu_{t_k}(dx) - \int_r^\infty dx\,  u(x,t_k)\Big| < \eps \Big\}
     \end{align}
        or, in other words,
	\begin{equation}
	\label{DD2.4}
\mathcal X_{T,\eps,n}(u) = \bigcap_{k=0}^{2^n} \Big\{
 u(r,t_k)dr \preccurlyeq \mu_{t_k}(dr) \preccurlyeq  u(r,t_k)dr \;\, \text{modulo $\eps$},
\Big\}
     \end{equation}

\smallskip

Recall that $P^{(N)}$ is the law of the true process $\und x(t)$ with initial distribution given by $N$ independent copies of a position variable distributed with probability $\rho_0(r) dr$.

	\begin{theorem}
\label{thmDD2.1}
There are a continuous probability density function $u(r,t)$, $t\ge 0$, $\alpha>0$ and $c$ so that:
	\begin{align}
	\nn
	&&\hskip-2.5cm
 P^{(N)}\big[\mathcal X_{T,\eps,n}(u)\big] \ge 1 -c 2^n N^{-\alpha}, \quad \text{ $T,\eps,n,N$ such that:}
\\&&  \eps \ge c\Big(\frac{e^T}{\sqrt T} 2^{3n/2}N^{-\alpha}+e^{2T}\frac {T^{N^{1/24}}}{N^{1/24}!}+ 2^{-n}T\Big)
 \label{DD2.3}
\end{align}

\end{theorem}

\smallskip

Thus  when $N\to \infty$ we can  take $T$ (the time window where we study the process) to $\infty$, the time grid $2^{-n}T$ (where we actually observe the process) as well as the closeness accuracy $\eps$ both to 0, provided that the last condition in \eqref{DD2.3} is satisfied.

We can characterize the limit density $u(r,t)$ as (1) the solution of a free boundary problem, if this has a classical solution, and (2) the unique separating element (in the sense of mass transport order) of ``lower and upper deterministic barriers'' denoted by $\rho^{\delta,\pm}(r,t)$.
We start from the former characterization.  Given a $C^1$ curve $\ga_t$, $t\ge 0$, we consider the evolution equation
	 \begin{eqnarray}
	\nn
&& \frac{\partial}{\partial t} \rho(r,t)= \frac 12 \frac{\partial^2}{\partial r^2}\rho(r,t)
+ \int dy \;\rho(y,t) p(y,r),\;\; r > \ga_t
\\  \label{DD2.5}
\\&& \rho(r,0)= \rho_0(r),\quad
\lim_{r\downarrow \ga_t}\rho(r,t)=0 \nn
   \end{eqnarray}
We will prove existence and uniqueness (see Theorem \ref{thmDD3.3}) and that  the solution of \eqref{DD2.5}
has a probabilistic interpretation similar to that in \eqref{1.3}.
The associated free boundary problem $\mathcal P$ 
is the following.

\medskip
$\mathcal P$: Let  $\rho_0(r)$ be as in the previous section, namely a smooth probability density with compact support; call
$\ga^*: =\inf\{r: \rho_0(r)>0\}$ 
	\begin{equation*}
\rho_0(\ga^*)=0,\qquad \lim_{r\downarrow \ga^*}\frac 12  \frac{d \rho_0(r)}{dr}=\int dy \;\rho_0(y,t) p(y,\ga^*)
	\end{equation*}
Find $\ga_t\in C^1$  and $\rho(r,t)$, $r \ge \ga_t$ in such a way that $\ga_0=\ga^*$,
$\rho(r,t)$ solves \eqref{DD2.5} and
  	 \begin{eqnarray}
  \label{DD2.6}
&&
\inf\{r\ge \ga_t: \rho(r,t)>0\} = \ga_t,\quad  \int_{\ga_t}^\infty dr\, \rho(r,t) =1 \quad \text{for all $t\ge 0$}
   \end{eqnarray}

\medskip

We will prove in Subsection \ref{subsec2.3} that:
	\begin{theorem}
\label{thmDD2.2}
If the free boundary problem $\mathcal P$ has a classical solution $(\ga_t, \rho(r,t)$) then  $\rho(r,t)$  coincides with the function $u(r,t)$ found in Theorem \ref{thmDD2.1}

\end{theorem}
In a recent paper, \cite{JM2},  Jimyeong Lee has proved local in time existence of classical solution of the free boundary problem $\mathcal P$ hence Theorem \ref{thmDD2.2} is not empty.
\medskip

The second characterization of the limit density $u(r,t)$ is actually the way
we prove its existence. In fact  $u(r,t)$ is obtained after introducing lower and upper barrier, denoted by $\rho^{\delta,\pm}(r,t)$, $\delta>0$. Their definition is given in Section \ref{subsec2.2}, their main property is:

	\begin{theorem}
\label{thmDD2.3}
For any $\delta$
positive  and any $t = k\delta$, $k$ positive integer:
	 \begin{equation}
	 \label{DD2.7}
\rho^{\delta,-}(r,t)dr \preccurlyeq \rho^{\delta,+}(r,t)dr
   \end{equation}
	 \begin{equation}
	 \label{DD2.8}
\rho^{\delta,-}(r,t)dr \preccurlyeq \rho^{\delta/2,-}(r,t)dr,\quad
\rho^{\delta/2,+}(r,t)dr \preccurlyeq \rho^{\delta,+}(r,t)dr
   \end{equation}
and finally, there is a unique function $u(r,t)$ such that
	 \begin{equation}
	 \label{DD2.9}
\rho^{\delta,-}(r,t)dr \preccurlyeq u(r,t)dr\preccurlyeq \rho^{\delta,+}(r,t)dr\quad \text{for all $\delta$ }
   \end{equation}
 The barriers are close in the $L^1$-norm: let $\delta =2^{-n}T$ then
    \begin{equation}
	 \label{DD2.9b}
\int dr\, \big|\rho^{\delta,-}(r,2^n \delta)- \rho^{\delta,+}(r,2^n \delta)| \le ce^T \delta
   \end{equation}

\end{theorem}

We give a detailed proof of \eqref{DD2.9b} in Subsection \ref{subsec2.2}
 and
 refer to the literature for the proof of the other statements in Theorem \ref{thmDD2.3} as very analogous to those in \cite{CDGPsurvey} and \cite{CDGP}
\cite{DF}, \cite{dfp} for similar models.

With the notation of  Theorem \ref{thmDD2.1}, in Section \ref{sec4} we will prove that
	 \begin{equation}
	 \label{DD2.10}
\rho^{\delta,-}(r,t_k)dr \preccurlyeq \frac 1N \pi_{\und x(t_k)}(dr)\preccurlyeq \rho^{\delta,+}(r,t_k)dr,\quad \text{modulo $\eps'$}
   \end{equation}
with probability $\ge 1-c 2^n  N^{-\alpha}$; in \eqref{DD2.10}  $\delta= 2^{-n}T$ and
$\eps' = c \frac{e^T}{\sqrt T} 2^{3n/2}N^{-\alpha}$.
The proof of \eqref{DD2.10} is in two steps. In the first step we introduce (see Subsection \ref{subsec2.1}) new processes $\und x^{\delta,\pm}$, called ``the stochastic upper and lower barriers'' and prove in Theorem \ref{thm1} that there are couplings so that for all $t_k$:
$$
 \pi_{\und x^{\delta,-}(t_k) } \preccurlyeq  \pi_{\und x(t_k)  }
 $$
with probability 1 as well as
$$
 \pi_{\und x(t_k)  } \preccurlyeq \pi_{\und x^{\delta,+}(t_k) }
 $$
In the second step we prove that
	 \begin{equation*}
\rho^{\delta,\pm}(r,t_k)dr \preccurlyeq \frac 1N \pi_{\und x^{\delta,\pm}(t_k)}(dr)\preccurlyeq \rho^{\delta,\pm}(r,t_k)dr,\quad \text{modulo $\eps'$}
   \end{equation*}
and thus get \eqref{DD2.10}.

We next use Theorem \ref{thmDD2.3}:  by \eqref{DD2.9b}
  	 \begin{equation*}
\rho^{\delta,+}(r,t)dr \preccurlyeq \rho^{\delta,-}(r,t)dr\quad \text{modulo $c \delta$}
   \end{equation*}
hence by \eqref{DD2.9}
 	 \begin{equation*}
\rho^{\delta,+}(r,t)dr \preccurlyeq u(r,t)dr\quad \text{modulo $c \delta$}
   \end{equation*}
Since an analogous argument holds for the lower barriers we get \eqref{DD2.3}  with $u(r,t)$ the function
 which separates upper and lower barriers;
the last term in \eqref{DD2.3} being the term $c\delta$ above.

To  prove Theorem \ref{thmDD2.2} we will show in Theorem \ref{thm3} that a classical solution of the free boundary problem $\mathcal P$ (see \eqref{DD2.6} ) is squeezed in between the lower and upper barriers and it thus coincides with their separating element $u(r,t)$.

\noindent
{\bf Outline of the paper.}
In Section \ref{sec:DD3} we will define the barriers and prove in Subsection \ref{subsec2.3} Theorem \ref{thmDD2.2} by exploiting a probabilistic representation of the solution of the free boundary problem which is related to one used in the definition of the barriers. In Section \ref{sec:DD4} we define the stochastic upper and lower barriers and prove  that they squeeze in between the true process (in the sense of mass transport order).
We will prove in Section \ref{sec4} Theorem \ref{thmDD2.1} using estimates whose proofs are postponed to Section \ref{sec:DDF.7}. In Section \ref{sec5} and in an Appendix we prove some more technical estimates.

\section{Probabilistic representations of deterministic evolutions}
\label{sec:DD3}
In this section we will first study \eqref{DD2.5} in the whole $\mathbb R$, namely the equation
 \begin{equation}
 \label{DD3.5}
\frac{\partial}{\partial t} \rho(y,t)=L^*\rho(\cdot,t)(y),\quad
 L^*f(y)=\frac 12 \frac{\partial^2f(y)}{\partial y^2}
+ \int dz \;f(z)p^*(y,z)
   \end{equation}
where $p^*(z,z')=p(z',z)$, referring to \eqref{DD3.5} as the {\em {free evolution}}.  We will then extend the analysis to  \eqref{DD2.5} itself, supposing that
$\ga_t$ is a  $C^1$ curve. The important point in both cases is a probabilistic representation of the solutions which will be often used in the sequel.

After that we will  define the ``lower and upper barriers''. Using the previous analysis and in particular the probabilistic representation of the solution of \eqref{DD2.5}
we will conclude the section by proving that a classical solution of the free boundary problem, when it exists, coincides with the separating element of Theorem \ref{thmDD2.3}.

\subsection{The free evolution}
\label{subseDD3.1}

As mentioned the solutions of \eqref{DD3.5} have a probabilistic representation.  We start from the latter and introduce the process $X_t$, $t\ge 0$, which is a brownian motion on $\mathbb R$ with jumps occurring independently at rate 1 with probability $p(x,y)$.  This is therefore a Markov process with generator
	\begin{equation}
	 \label{DDD3.1}
\mathcal Lf(x)=\frac 12 \frac{d^2 f(x)}{d x^2}+\int p(x,y)[ f(y)-f(x)]dy
	\end{equation}
$f\in C^2$. $\mathcal L$ defines a semi-group on $ C(\mathbb R,\mathbb R_+)\cap L^\infty$ denoted by  $S_t=e^{\mathcal L t}$.  If $P_x$ is the law of the process $X_t$ starting from $X_0=x$ and $E_x$ its expectation,  then if $f\in  C(\mathbb R,\mathbb R_+)\cap L^\infty$
	\begin{equation}
	 \label{DDD3.2}
E_x\big[ f(X_t\big]= e^{\mathcal L t}f(x)
	\end{equation}
We will also consider the adjoint process $X^*_t$ with generator
	\begin{equation}
	 \label{DDD3.3}
\mathcal L^*f(x)=\frac 12 \frac{\partial^2 f(x)}{\partial x^2}+\int p^*(x,y)[ f(y)-f(x)]dy, \quad p^*(x,y)=p(y,x)
	\end{equation}
and the corresponding semi-group $S^*_t=e^{\mathcal L^* t}$ on $ C(\mathbb R,\mathbb R_+)\cap L^\infty$.

  \begin{theorem}
\label{thmDD3.1}
The Markov semi-groups $e^{\mathcal L t}$ and $e^{\mathcal L^* t}$ have a kernel
denoted respectively by $e^{\mathcal L t}(x,y)$ and $e^{\mathcal L^* t}(x,y)$ which are $C^\infty$ for $t>0$ and differentiable in $t$. Call $T_t(x,y):=e^te^{\mathcal L t}(x,y)=e^te^{\mathcal L^* t}(y,x)=:T^*_t(y,x)$, then $T^*_t(y,x)$ is the Green function for \eqref{DD3.5}, namely for any $x$, $T^*_t(y,x)$, as a function of $(y,t)$, solves  \eqref{DD3.5} with initial datum $\delta_x(r)$.  In particular the solution $\rho$
of  \eqref{DD3.5} with initial datum $\rho_0$ is
	\begin{equation}
	 \label{DDD3.3.0}
\rho(x,t) = \int dy\,  T^*_t (x,y) \rho_0(y)=:T^*_t\rho_0(y),\quad \|\rho(\cdot,t)\|_\infty \le e^t\, \|\rho_0\|_\infty 
	\end{equation}
  \end{theorem}

The proof of Theorem \ref{thmDD3.1} is elementary but we will give for completeness some details below.  Before that we set some notation and remark some consequences of the Theorem.

We  denote by $P_{x,s}$, $s \ge 0$, the law of the process $\{X_t, t \ge s\}$ with $X_s=x$ and write
$P_{x,s}\big(X_t\in dy\big)$ for its law. When $s=0$ we simply write $P_x\big(X_t\in dy\big)$. From Theorem \ref{thmDD3.1} it follows that this law is absolutely continuous with respect to the Lebesgue measure:
 \begin{equation}
	 \label{DD3.6}
e^t P_x\big(X_t\in dy\big) =T_t(x,y)dy,\qquad e^{t-s} P_{x,s}\big(X_t\in dy\big) = T_{t-s}(x,y)dy
	\end{equation}
The proof of \eqref{DDD3.3.0} in Theorem \ref{thmDD3.1} follows from the others statements. In fact  for $u$ and $f$  in $C(\mathbb R,\mathbb R_+)\cap L^\infty$ and $\rho(r,t)$ solution of \eqref{DD3.5} with initial datum $\rho_0$ we have
  	\begin{eqnarray}
	\nn
\int dx\, \rho(x,t) f(x) &=& e^t\int dx\,  \int dy\,  \rho_0(x) e^{\mathcal L t}(x,y) f(y) \\&=&e^t \int dy \,  \int dx\, f(y) e^{\mathcal L^* t}(y,x) \rho_0(x)
 \label{DDD3.4}
	\end{eqnarray}
which  by the arbitrariness of $f$ gives \eqref{DDD3.3.0}. \eqref{DDD3.4} expresses the well known fact that  we can compute the expectation of a function $f$ at time $t$ either by evolving the initial measure $\rho_0(r)dr$ till time $t$ or by evolving $f$ and then computing the expectation with respect to the initial measure: the former is obtained using the semi-group $e^{\mathcal L^* t}$ the latter with the semi-group $e^{\mathcal L t}$.

By \eqref{DDD3.4} the solution of \eqref{DD3.5} is given by the semigroup $e^te^{\mathcal L^* t}=: e^{L^* t}$, where $L^*= \mathcal L^*+1$ is made explicit in \eqref{DD3.5}.
We also denote by $L= \mathcal L+1$ the adjoint of $L^*$.

In the next lemma we establish the existence and uniqueness of the Green function for  \eqref{DD3.5}.

\begin{lemma}
	\label{lemDD3.1}
Given any $x\in\mathbb R$ there is a unique solution $v(y,t)$ of \eqref{DD3.5} with initial condition $v(y,0)=\delta_x(y)$. Supposing $p(0,x)$ in $C^\infty$ also $v\in C^\infty(\mathbb R,(0,\infty))$. Furthermore
there are $c$ and $c'$
so that
	 \begin{equation}
	 \label{DD3.2}
\sup_r |\frac{\partial}{\partial r} v(r,t)| \le c\,\, \frac {e^t}{ t},\quad
 |\frac{\partial}{\partial t} v(r,t)| \le c'\, \Big(|r-x|^{-2}+e^t\Big)
   \end{equation}
   	 \begin{equation}
	 \label{DD3.222}
\int dr\, |\frac{\partial}{\partial r} v(r,t)| \le c ( \frac 1{\sqrt t} + e^t)
   \end{equation}
 	\end{lemma}

\noindent{\bf Proof.} We start by solving the integral version of \eqref{DD3.5}.
By an abuse of notation we   call  $G_t$ and $p^*$ the integral operators with kernels $G_t(x,y)$ and $p^*(x,y)$ respectively, $G_t(x,y)=\frac 1{\sqrt{2\pi t}} \exp\{-\frac{(x-y)^2}{2t}\}$.  Then
\begin{equation*}
v(y,t)=G_t\delta_x(y)+\int_0^t ds\,G_s p^* v(\cdot,t-s)(y)
	\end{equation*}
explicitly:
\begin{equation}
	\label{aDD3.1}
v(y,t)=G_t(y,x)+\int_0^t ds  \int dz \,G_s(y,z) \int dz' \, p^*(z,z')v(z',t-s),\qquad x,y\in\mathbb R
	\end{equation}
Iterating \eqref{aDD3.1} we get
	\begin{eqnarray}
	v(y,t)= G_t\delta_x(y)+\sum_{n=1}^\infty \int_{\mathbf R^n_+}ds_1..ds_n \mathbf 1_{s_1\le s_2\le..\le s_n\le t}\Big(G_{s_1}p^*G_{s_2}p^*\cdots p^*G_{t-s_n}\delta_x\Big)(y)\nn\\ \label{DD3.3}
	\end{eqnarray}
  After bounding by $c:=\sup_x p(0,x)$ the last factor $p^*$ on the right hand site of \eqref{DD3.3} we get
	\begin{eqnarray}
	 \label{DD3.3.0}
	|v(y,t)| \le  G_t(y,x)+c(e^t-1)
	\end{eqnarray}
Thus $v$ as given by \eqref{DD3.3} is well defined and it is the unique solution of \eqref{aDD3.1}.
Furthermore by differentiating \eqref{DD3.3} $k$ times with respect to $y$ we get
\begin{eqnarray*}
	\|  v^{(k)} \|_\infty\le  \|G_t^{(k)}\|_\infty+e^t \|p\|_\infty\|p^{(k)}\|_\infty 2\xi
		\end{eqnarray*}
with $\xi$ the range of $p(0,\cdot)$.
We
differentiate  \eqref{aDD3.1}
with respect to $t$, getting
	\begin{eqnarray*}
\frac{\partial}{\partial t}	v(y,t) &=& \frac{\partial}{\partial t} G_t(y,x)+
\int dz'\,  p^*(y,z') v(z',t) \\ &+& \int_0^t \int dz\,\int dz'\,  \frac 12	G_{t-s}(y,z)\frac{\partial^2}{\partial z^2}p^*(z,z') v(z', s)
	\end{eqnarray*}
By iteration we prove that all the derivatives with respect to $t$ of $v$ exist.
By \eqref{DD3.3.0} and recalling that $p^*(z,z')$ is smooth with compact support we get that for suitable constants $c_1$ and $c_2$:
	\begin{eqnarray*}
|\frac{\partial}{\partial t}	v(y,t)| &\le& |\frac{\partial}{\partial t} G_t(y,x)|+
c_1 e^t  + c_2(t^{1/2}+ e^t)
	\end{eqnarray*}
hence \eqref{DD3.2}.

Proof of \eqref{DD3.222}.  By \eqref{DD3.3} and denoting by $v_y(y,t)$ the derivative with respect to $y$:
	\begin{eqnarray*}
&&	\int dy \,|v_y(y,t)| \le \frac 1{\sqrt{2\pi}}\int dy\, \frac {|x-y|}{\sqrt{t^{3/2}}} e^{-\frac{(x-y)^2}{2t}} \\ && \hskip2cm+ \sum_{n=1}^\infty \int_{\mathbf R^n_+}ds_1..ds_n \mathbf 1_{s_1\le s_2\le..\le s_n\le t}\int dy\Big(G_{s_1}[2\xi \|p'\|_\infty] G_{s_2}p^*\cdots p^*G_{t-s_n}\delta_x\Big)(y)
	\end{eqnarray*}
The first term is bounded by $c \frac 1{\sqrt t}$, the second one proportionally to $e^t$.

It is easy to see that any solution of the differential equation \eqref{DD3.5} is also a solution of the integral equation  \eqref{aDD3.1}, thus the lemma is proved.
\qed

Call $K_t(x,y):=v(y,t)$, $v$ as in Lemma \ref{lemDD3.1}, then
  $$
 \rho(x,t):=\int dy\, K_t(x,y) \rho_0(y)
 $$
is the unique solution of \eqref{DD3.5} with initial datum $u$. Denoting by $K_t$ the integral operator with kernel $K_t(x,y)$, we have $\rho(r,t)= K_t \rho_0(r)$ and by
  \eqref{DD3.5}, $\frac{d}{dt} K_t \rho_0 = L^*K_t \rho_0$,
hence $K_t = e^{L^*t}=e^t e^{\mathcal L^* t}$ and
Theorem \ref{thmDD3.1} is proved.

\subsection{The  evolution in semi-infinite domains}
\label{subseDD3.2}
	
There is also a probabilistic representation for the solution of \eqref{DD2.5} for a given  $C^1$- curve $\Ga=\{\ga_t, t\ge 0\}$.  Let $X_t$ be the process defined in the previous subsection
and call  $\tau_s=\inf\{t>s: X_t\le \ga_t\}$. Given $s\ge 0$ and $X_s>\ga_s$ we define for $t>s$
	\begin{eqnarray}
X^{\Ga}_{t;s}=\begin{cases}X_t &\text{ if }\tau_s>t
\\ -\infty &\text{ otherwise }
	\end{cases}
\end{eqnarray}
For $s\ge 0$ and $x>\ga_s$  we call
$P_{x,s}\big(X^{\Ga}_t\in dy\big)$ the law  at time $t>s$ of $X^{\Ga}_t$ on $\mathbb R$ with $X^{\Ga}_s=x$   and claim that it has a density with respect to the Lebesgue measure, denoted by  $\alpha_{x,s}(y,t)$ .  This follows because by \eqref{DD3.6}
	\begin{equation}
     \label{n3.12}
P_{x,s}\big(X^{\Ga}_t\in dy\big)=P_{x,s}\big(X_t\in dy;\tau_s>t\big)\le P_{x,s}\big(X_t\in dy\big) = e^{-(t-s)} T_{t-s}(x,y)dy
	\end{equation}
We will simply write $\alpha_x(y,t)$ for $\alpha_{x,0}(y,t)$.

\medskip
The main result of this Subsection is the following theorem.

 \begin{theorem}
\label{thmDD3.3}
For all $x\in\mathbb R$ and all $t>s_0\ge 0$ $\alpha_{x,s_0}(\cdot,t)\in C^\infty((\ga_t,\infty),\mathbb R_+)$ and it  is differentiable in $t$.  Furthermore  $v_{x,s_0}(y,t):=e^{t-s_0}\alpha_{x,s_0}(y,t)$
 satisfies \eqref{DD2.5} and  $v_{x,s_0}(y,s_0)=\delta_x(y)$.
  \end{theorem}

We split the proof of Theorem \ref{thmDD3.3} in three lemmas: in Lemma \ref{lemDD3.2} we  prove that $\alpha_{x,s_0}(y,t)$ is smooth, in Lemma \ref{lemDD3.3} we prove that $ \alpha_{x,s_0}$  is 0 at the boundary and finally in Lemma \ref{lemDD3.4} we compute the  derivative with respect to $t$ of $\alpha_{x,s_0}(y,t)$. All that proves Theorem \ref{thmDD3.3}.

\begin{lemma}
\label{lemDD3.2}
For all $t>s_0\ge 0$ $\alpha_{x,s_0}(\cdot,t)\in C^\infty((\ga_t,\infty),\mathbb R_+)$, it  is differentiable in $t$ and  for all $x>\ga_0$
	\begin{equation}
	\label{DD3.8}
\alpha_{x,s_0}(y,t)=e^{-(t-s_0)} T_{t-s_0}(x,y)-\int_{s_0}^t \int_{z\le \ga_s}q_{x,s_0}(ds dz)e^{-(t-s)}T_{t-s}(z,y),\quad y>\ga_t
	\end{equation}
where $q_{x,s_0}(ds dz)=P_{x,s_0}(\tau_{s_0}\in ds, X_s\in dz)$.  Moreover 	\begin{equation}
	\label{DD3.8.1}
\sup_{x> \ga_{s_0},y>\ga_t} \alpha_{x,s_0}(y,t)< \infty
	\end{equation}
\end{lemma}

\noindent{\bf Proof.} For notational simplicity we take $s_0=0$.
Call $A_x(y,t)$ the right hand side of \eqref{DD3.8}. By Lemma \ref{lemDD3.1} $A_x(y,t)$ is a smooth function of $y$ in $\{y>\ga_t\}$ and it is differentiable in $t$. Thus we only need to prove that $A_x(y,t)=\alpha_{x}(y,t)$. For all $f$ with support on $y>\ga_t$
	\begin{equation*}
\mathbb E_x(f(X^\Ga_t))=E_x(f(X_t)\mathbf 1_{\tau_0>t})=E_x(f(X_t))-E_x(f(X_t)\mathbf 1_{\tau_0\le t})
	\end{equation*}
By conditioning
	\begin{equation*}
E_x(f(X_t)\mathbf 1_{\tau_0\le t}) = \int_0^t \int_{z\le \ga_s}q_x(ds dz)\int E_{z} [ f(X_{t-s})]
	\end{equation*}
Then using \eqref{DD3.6}
	\begin{equation*}
\int dy\,\alpha_{x}(y,t) f(y)= \int dy\,A_{x}(y,t) f(y)
	\end{equation*}
and \eqref{DD3.8} follows by the arbitrariness of $f$. \eqref{DD3.8.1} follows from \eqref{n3.12} after using
 Lemma \ref{lemDD3.1} to bound
 $ T_{t-s_0}(x,y)$.\qed

\begin{lemma}
\label{lemDD3.3}
For all $t>s_0\ge 0$ and $x> \ga_{s_0}$
\begin{equation}
	\label{DD3.9}
\lim_{z\to \ga_t} \alpha_{x,s_0}(z,t) =0
	\end{equation}
\end{lemma}

\noindent{\bf Proof.}  For notational simplicity we take $s_0=0$.
Let $0<s<t$ then, by conditioning
$
\alpha_{x}(z,t)= \int_{\ga_s}^\infty dy\,
\alpha_{x}(y,s) \alpha_{y,s}(z,t)
$.
Using the reverse process we get
 \begin{equation*}
	\label{DD3.9.2}
\alpha_{x}(z,t)= E^*_{z}\big[\alpha_x(X^* _{t-s},s); X^*_{s'} > \ga_{t-s'}, s'\in[0,t-s]\big]
	\end{equation*}
where $E^*$ denotes expectation with respect to the process $X^*_t$ of a brownian motion which also jumps (from $x$ to $y$) at rate 1 with probability $p^*(x,y)dy$.  By \eqref{DD3.8.1} there is $c$ so that the last term is bounded by
\begin{equation*}
	\label{DD3.9.3}
\alpha_{x}(z,t) \le c  P^*_{z}\big[ X^*_{s'} > \ga_{t-s'}, s'\in[0,t-s]\big]
	\end{equation*}
Since $\ga_t$ is $C^1$ the right hand side vanishes as $z \to \ga_t$.  \qed

\begin{lemma}
\label{lemDD3.4}
Let $t>s_0\ge 0$.
For all $y>\ga_t$
\begin{eqnarray}
	\label{DD3.10}
\frac{\partial \alpha_{x,s_0}(y,t)}{\partial t}&=& \mathcal L^*  \alpha_{x,s_0}(y,t)
\\ &=&
\frac 12 \frac{\partial^2 \alpha_{x,s_0}(y,t)}{\partial y^2}+\int_{\ga_t}^\infty dz\, p(z,y)\{\alpha_{x,s_0}(z,t) - \alpha_{x,s_0}(y,t)\} \nn
	\end{eqnarray}

\end{lemma}

\noindent{\bf Proof.}
For notational simplicity we take $s_0=0$.
Let $b>a>\ga_t$ and let $\varphi$ be a smooth function with support in $[a,b]$. Then
	\begin{eqnarray*}
	\nn
\int dy\, \phi(y) \frac{\partial \alpha_{x}(y,t)}{\partial t} = \lim_{\eps\to 0} \int dz\,\, \alpha_x(z,t)
\frac 1 \eps \Big( \int dy\,\alpha_{z,t}(y,t+\eps)\varphi(y)-\varphi(z)\Big)
	\end{eqnarray*}
By \eqref{DD3.8}
	\begin{eqnarray*}
	\nn
\int dy\, \alpha_{z,t}(y,t+\eps)\varphi(y)-\varphi(z) =\{ e^{\mathcal L \eps}\varphi(z)-\varphi(z)\}
- \int_0^\eps \int q_{z,t}(ds dz')e^{\mathcal L (\eps-s)}\varphi(z')
	\end{eqnarray*}
We are going to prove that
	\begin{eqnarray}
	\label{3.100}
 \lim_{\eps\to 0}\int dz\,\, \alpha_x(z,t)
\frac 1 \eps \int_0^\eps \int q_{z,t}(ds dz')e^{\mathcal L (\eps-s)}\varphi(z')=0
	\end{eqnarray}
To   this end we split the integral over $z$ into $z> \ga_t+\delta$ and $z \in [\ga_t,\ga_t+\delta]$.  When  $z> \ga_t+\delta$ the main contribution to $q_{z,t}(ds dz')$
comes from a jump which happens with probability of order $\eps$ (the probability that a brownian reaches the curve, which is $C^1$, in a time $\le \eps$ starting from $z> \ga_t+\delta$ is of higher order).  For the same reason $e^{\mathcal L (\eps-s)}\varphi(z')$, $z' \le \ga_s$ has order $\eps$ hence the contribution to \eqref{3.100} from $z> \ga_t+\delta$ vanishes.  Using again that  $e^{\mathcal L (\eps-s)}\varphi(z')$, $z' \le \ga_s$ has order $\eps$ since   $\alpha_x(z,t)$ is bounded then  the contribution of the integral over $z$ from the interval $ [\ga_t,\ga_t+\delta]$ is bounded by $c \delta \eps$.  By the arbitrariness of $\delta$
this proves \eqref{3.100}.

We have thus shown that
	\begin{eqnarray*}
	\nn
\int dy\, \phi(y) \frac{\partial \alpha_{x}(y,t)}{\partial t} =  \int dz\,\, \alpha_x(z,t)
\mathcal L \varphi(z)
	\end{eqnarray*}
hence \eqref{DD3.10} by the arbitrariness of $\varphi$. \qed

\subsection{The deterministc barriers}
		\label{subsec2.2}
In this Subsection we  define the  {\it {deterministic barriers}} $\rho^{\delta,\pm}(\cdot,t)$, $\delta>0$, which appear in Theorem \ref{thmDD2.3}. We first define $\rho^{\delta,\pm}(x,k\delta)$ for any positive integer $k$.
To this end we first introduce the {\it {cut operators}} $C^\pm_\delta$ as follows.
 Calling $\mathcal M_a$, $a>0$, the set of $v\in L^1(\mathbb R,\mathbb R_+)$ such that $\int dx\, v(x)=a$ we define $C^-_\delta: \mathcal M_1 \to \mathcal M_{
 e^{-\delta}}$ as the operator which cuts on the left a mass $1-e^{-\delta}$. Analogously $C^+_\delta: \mathcal M_{
 e^{\delta}} \to \mathcal M_{1}$ cuts on the left a mass $e^{\delta}-1$.  Thus
 \begin{eqnarray}
	 \label{2.5a}
&&C^+_\delta
v(x)=\mathbf 1_{x\ge V^+}\; v(x),\qquad V^+\text{ such that } \int_{V^+}^\infty dx\, v(x)=1\nn \\&&
\\&& C^-_\delta v(x)=\mathbf 1_{x\ge V^-}\; v(x),\qquad V^-\text{ such that } \int_{V^-}^{\infty} dx\,v(x)=e^{-\delta}\nn
   \end{eqnarray}
Let $T^*_t$ as in \eqref{DDD3.3.0}, then $T^*_\delta C^-_\delta: \mathcal M_1 \to \mathcal M_1$ and $C^+_\delta T^*_\delta : \mathcal M_1 \to \mathcal M_1$ and 
for  any $v\in \mathcal M_1$ and any integer $k>0$ we then set
     \begin{equation}
    \label{2.6a}
\rho^{\delta,-}(x,k\delta) =(T^*_\delta C^-_\delta)^kv(x) ,\qquad \rho^{\delta,+}(x,k\delta) =(C^+_\delta T^*_\delta)^kv(x)
\end{equation}
We finally define
	 \begin{equation}
    \label{2.6b}
\rho^{\delta,\pm}(x,t) =T^*_{t- k\delta}\rho^{\delta,\pm}(x,k \delta) \qquad  \text{ for }\,k\delta< t<(k+1)\delta
\end{equation}

{\bf Proof of \eqref{DD2.9b}.}    Fix $\delta>0$ and shorthand
$T^*=T^*_\delta$, $C^{\pm}=C^{\pm}_\delta$.  We want to prove that for any  $n>0$ and any $v\in \mathcal M_1$,
	 \begin{equation}
    \label{4.1000.1}
\| (C^+T^*)^nv - (T^*C^-)^n v\|_{L^1} \le e^{n\delta}(e^\delta-e^{-\delta})
+2 e^\delta(e^\delta-1)
\end{equation}
For any  $w\in \mathcal M_1$
	 \begin{equation}
    \label{4.1000.2}
e^\delta C^- w= C^+ e^{\delta}w
    \end{equation}
and since $ e^{\delta}T^*=T^* e^{\delta}$ we have
	 \begin{eqnarray}
    \label{4.1000.3}
&&(T^*C^-)^n v=  e^{-\delta}\phi,\quad \phi:= T^* \Big((C^+T^*)^{n-1}\Big)  e^\delta C^-v
\in \mathcal M_{e^\delta}\nn\\&&
\\&&
(C^+T^*)^n v =  C^+ \psi, \quad \psi:=T^* \Big((C^+T^*)^{n-1}\Big) v \in \mathcal M_{e^\delta}\nn
\end{eqnarray}
By \eqref{2.5a}
	 \begin{equation}
    \label{4.1000.5}
\|e^\delta C^-v - v\|_{L^1} \le 1-e^{-\delta}+ \|(e^\delta-1)v \|_{L^1} =e^\delta-e^{-\delta}
\end{equation}
so that
	 \begin{equation}
    \label{4.1000.6}
\|\psi-\phi \|_{L^1} \le e^{\delta n}\Big(e^\delta-e^{-\delta}\Big)
\end{equation}
By \eqref{4.1000.3}
	 \begin{equation}
    \label{4.1000.7}
\| (C^+T^*)^nv - (T^*C^-)^n v\|_{L^1} \le \|e^{-\delta} \phi - C^+ \psi\|_{L^1}
\end{equation}
and finally, recalling that $\psi$ and $\phi$ are in $\mathcal M_{e^\delta}$,
	 \begin{equation*}
\|e^{-\delta} \phi - C^+ \psi\|_{L^1} \le \|\psi-\phi \|_{L^1} +(1-e^{-\delta} ) e^\delta
+( e^\delta-1) e^\delta
\end{equation*}

\qed

\subsection{A priori bounds}
\label{subsec2.2b}

We will  prove below that if the initial datum $\rho_0$ is bounded then both the solution of \eqref{DD3.5} and the barriers are bounded. To prove such a result we will use the probabilistic representation of the evolution.

\begin{lemma}
\label{lemmaDD3.2}
Let $\rho(x,t)$ be the solution of \eqref{DD3.5} with initial condition $\rho_0$.
For any  $b>0$  there is $c_{b}$ so that for any $t>0$
	 \begin{equation}
    \label{DD5.12}
\int_{|x| \ge N^b}dx\, \rho(x,t) \le c_b e^t\frac{t^{N^{b/2}}}{N^{b/2} !}
 \end{equation}
\end{lemma}

\noindent
{\bf Proof.}  Let
$[-A,A]$ be the support of  $\rho_0$, then by \eqref{DDD3.3.0}
	 \begin{equation*}
\int_{|x| \ge N^b}dx\, \rho(x,t) = \int_{-A}^A dx\,\rho_0(x)\int_{|r'| \ge N^b}dy\,
e^te^{\mathcal Lt}(x,y)
     \end{equation*}
Let $[-\xi,\xi]$ be the support of $p(0,y-x)$
then the above integral is bounded by the sum of (1) the probability that the number $n(t)$ of jumps
is  $
   n(t) \ge \frac{N^b - A}{2\xi}\ge N^{b/2}
   $
and (2) a brownian motion $B(t)$ starting from 0 is   $|B(t)|\ge \frac{N^b - A}{2}$.  \eqref{DD5.12} then follows.  \qed

The fact  that the barriers are bounded is a consequence of the following Lemma.

\begin{lemma}
\label{lemmaDD3.7}
Let $\rho(x,t)$ be the solution of \eqref{DD3.5} with initial condition $\rho_0$.
Then the deterministic upper and lower barriers with same initial condition are bounded by $\rho(r,t)$:
	 \begin{equation}
    \label{DD3.24}
 \rho^{\delta,\pm}(x,t)\le \rho(x,t),\qquad \forall x\in\mathbb R, \quad t=k\delta
 \end{equation}
	\end{lemma}

\noindent
{\bf Proof.}  Observe that for all functions $f$ and $g$ we have that $T^*_t(f+g)(x)\ge T^*_t f(x)$ for all $x\in\mathbb R$. Also for all positive function $f$ we have that $C^\pm f\le f$. Thus \eqref{DD3.24} obviously holds.\qed

\subsection{ Proof of Theorem \ref{thmDD2.2} }
\label{subsec2.3}

 In this subsection we establish a relation between the  free boundary problem and the deterministic  barriers and
 prove Theorem \ref{thmDD2.2}  namely that the classical solution of the free boundary problem, when it exists,  is the separating element of the deterministic  barriers. To prove this it is enough to show that the classical solution is squeezed in between the lower and the upper barriers which is done in Theorem \ref{thm3} below.
The proof of the theorem is similar to others for analogous models, they all exploit the representation of solutions of the heat equation in terms of Brownian motions, in particular that the hitting probability of a Brownian at a curve $\ga_t$ has a density with respect to Lebesgue, a property which is well known for $C^1$ curves but which extends to Holder curves with parameter $>1/2$, see \cite{JM}.

Assume $(\rho(\cdot, t), \ga_t)_{t\in [ 0,T]}$ is a classical solution of the free boundary problem $\mathcal P$ (see   \eqref{DD2.6} recalling that the initial datum is a probability density $\rho_0\in   L^\infty(\mathbb R,\mathbb R_+)$.  Then

	\begin{theorem}
\label{thm3}
Consider the deterministic barriers with initial datum $\rho_0$, then for all  $\delta>0$ and all positive integer $k$ we have
	\begin{align}
	\label{2.9}
\rho^{\delta,-}(\cdot,k\delta)\preccurlyeq \rho(\cdot,k\delta)\preccurlyeq\rho^{\delta,+}(\cdot,k\delta)\end{align}

\end{theorem}

\noindent {\bf Proof.}
We use the probabilistic representation of Section \ref{subseDD3.2}, thus
  for any  $r\in \mathbb R$
	 \begin{equation}
	 \label{6.2}
\int_r^\infty \rho(x,t)dx=e^t\int \rho_0(x)P_x\big(X^{\Ga}_t>r\big)  dx =
e^t\int \rho_0(x)P_x\big(X_t>r;\tau>t\big)  dx
\end{equation}
Recall that by \eqref{DD2.6} the left hand side for $r=-\infty$ is equal to 1. Observe that this implies that the distribution of $\tau$ is exponential of parameter 1.
\medskip

\noindent{\bf Upper bound $k=1$.}
Recall that $\rho^{\delta,+}(x,\delta)=C^+T^*_\delta \rho_0(x)$. We  prove below that for any $r\in\mathbb R$
	\begin{equation}
	\label{6.3}
\int_r^\infty	\rho^{\delta,+}(x,\delta)dx\ge\int \rho_0(x) e^\delta P_x\big(X^{\Ga}_\delta \ge r\big)  dx
	\end{equation}

\noindent{\bf Proof}. Let $V^{+}$ be such that
$\rho^{\delta,+}(x,\delta)=T^*_\delta \rho_0(x)\mathbf 1_{x\ge V^{+}}$.
If $r<V^{+}$, then
\begin{align}
			 \nn
\int_r^\infty \rho^{\delta,+}(x,\delta)dx = \int_{V^{+}}^\infty T^*_\delta \rho_0(x)dx=1	\ge \int \rho_0(x) e^\delta P_x\big(X^{\Ga}_\delta\ge r\big) dx			\end{align}
If $r>V^{+}$
 		\begin{align}
			 \nn
\int_r^\infty \rho^{\delta,+}(x,\delta)dx &=e^\delta \int \rho_0(x) P_x\big(X_\delta>r\big)  dx           \ge e^\delta \int_{L_0}^\infty \rho_0(x) P_x\big(X_\delta>r;\tau>\delta\big)  dx
			\label{m.2.8}
			\end{align}
that by  \eqref{6.2} concludes the proof. \qed

\noindent{\bf Upper bound for all $k$.} Since
$
 \rho(\cdot,\delta)\preccurlyeq\rho^{\delta,+}(\cdot,\delta)$
and since $C_\delta^+T^*_\delta$ preserves  order (see Chapter 5 of \cite{CDGPsurvey})
\begin{align}
	\nn
C_\delta^+T^*_\delta  \rho(\cdot,\delta) \preccurlyeq  C_\delta^+T^*_\delta \rho^{\delta,+}(\cdot,\delta)
=\rho^{\delta,+}(\cdot,2\delta)\end{align}
Then  $
 \rho(\cdot,2\delta)\preccurlyeq\rho^{\delta,+}(\cdot,2\delta)$ is proved as in the case $k=1$. By iteration the argument extends to all $k$.

\noindent{\bf Lower bound. } We only prove it in the case $k=1$, the extension to all $k$ being similar to the one in the previous case.

We split $\rho_0=f+g$, $f(x)=\rho_0(x)\mathbf 1_{x \ge V^-}$, $g=\rho_0\mathbf 1_{x< V^-}$ with $V^-$ such that
	\begin{equation}
\label{4.500}
	\int f(x)dx=e^{-\delta},
	\end{equation}
Then $\rho^{\delta,-}(x,\delta)=T^*_\delta f(x)$ thus
	\begin{equation*}
	\int_r^\infty  \rho^{\delta,-}(x,\delta)dx= \int f(x) e^\delta P_x(X_\delta>r)
	\end{equation*}
Instead
 \begin{eqnarray*}
&&\hskip-.5cm
\int_r^\infty \rho(x,\delta)dx=e^\delta\int f(x)P_x\big(X^{\Ga}_\delta>r\big)  dx+ e^\delta\int g(x)P_x\big(X^{\Ga}_\delta>r\big)  dx\\&&= \int f(x) e^\delta P_x(X_\delta>r)
- e^\delta\int f(x)P_x\big(X_\delta>r:\tau\le \delta \big)  dx
e^\delta\int g(x)P_x\big(X^{\Ga}_\delta>r\big)  dx
\end{eqnarray*}
We thus need to prove that
\begin{eqnarray*}
 e^\delta\int g(x)P_x\big(X_\delta>r;\tau>\delta\big)  dx\ge e^\delta\int f(x)P_x\big(X_\delta>r;\tau\le \delta \big)  dx
\end{eqnarray*}
that we rewrite  as follows
\begin{eqnarray}
	 \label{6.4b}
 \int g(x)P_x(\tau>\delta) P_x\big(X_\delta>r|\tau>\delta\big)  \ge\int f(x)P_x(\tau\le \delta)P_{x} \big(X_\delta>r |\tau\le\delta\big)
\end{eqnarray}
First we prove that the two measures $g(x)P_x(\tau>\delta)dx$ and $f(x)P_x(\tau\le\delta)dx$ have same mass. The difference of the two masses is
\begin{eqnarray*}
\hskip-.4cm \int g(x)P_x\big(\tau>\delta\big)  dx-\int f(x)[1-P_x\big(\tau> \delta \big)]  dx 
=\int\rho_0(x)P_x\big(\tau>\delta\big)  dx-\int f(x)dx=0
\end{eqnarray*}
because from \eqref{6.2} we get $e^\delta\int\rho_0(x)P_x(\tau>\delta)=1$ while $\int f(x)dx=e^{-\delta}$ by \eqref{4.500}.
 We rewrite \eqref{6.4b} as
\begin{eqnarray*}
	 \label{6.4}
 \int \mu(dx)P_x\big(X_\delta>r|\tau>\delta\big)  \ge\int \la(dx,dz,ds)P_{z,s} \big(X_\delta>r \big)
\end{eqnarray*}
where $\mu(dx)=g(x)P_x(\tau>\delta)dx$ and
$
\la(dx,dz,ds)= f(x)dxP_x(\tau\in ds, X_s\in dz)
$ namely the probability that the process hits the region $\{x\le \ga_t\}$ at the point $dz$ at time $ds$. Since $\mu$ and $\la$ have the same mass \eqref{6.4b} follows from
	$$
P_{z,s} \big(X_\delta>r \big)  \le P_x\big(X_\delta>r|\tau>\delta\big)  ,\qquad \text{ for all }s\in[0,\delta),  z\le \ga_s\text{ and }x\ge \ga_s$$
which can be proved  as in Section 10 of \cite{CDGPsurvey}. We omit the details.

\section{Stochastic barriers}
\label{sec:DD4}

\subsection{Definition of the stochastic barriers}
\label{subsec2.1}

For each positive real number $\delta$ we define two processes $\und x^{\delta,+}(t)$ and $\und x^{\delta,-}(t)$, $t\ge 0$, called respectively \emph{upper and lower stochastic barriers}.
We are going to define inductively $\und x^{\delta,\pm}(t)$.
We thus suppose to have defined $\und x^{\delta,\pm}(t)$ for $t\le t_{k-1}^+$  and want to define $\und x^{\delta,\pm}(t)$ for $t\le t_k^+$, $t_k=k\delta$.

\begin{itemize}
\item \emph{The upper stochastic barrier}. We set $\und x^{\delta,+}(0^+)=
\und x(0)$ and suppose inductively that we have defined the process till time
$t_{k-1}^+$, $k\ge 1$.
For $t \in [t_{k-1}^+,t_k^-]$ the process $\und x^{\delta,+}(t)$ is defined as the basic process $\und y(t)$ starting from $\und x^{\delta,+}(t_{k-1}^+)$, namely  the particles evolve as independent Brownian motions with non local branching.  Calling $\und x^{\delta,+}(t_k^-)$ the final positions of these particles (i.e.\ at time $t_k^-$), we then define
  \begin{align}
    \und x^{\delta,+}(t_k^+):= N \hbox{ rigthmost particles of }\und x^{\delta,+}(t_k^-).
  \end{align}

\item \emph{The lower  stochastic barrier}. The definition is again by induction and it depends on a number  $\alpha_0\in (1/2,1)$. Call
	\begin{equation}
	\label{2.2}
M_\delta:=(1-e^{-\delta})N + N^{\alpha_0}
	\end{equation}
Initially we set $\und x^{\delta,-}(0^+)$ as the configuration obtained from $\und x(0)$  by deleting the leftmost  $M_\delta$  particles and then define
$\und x^{\delta,-}(t)$, $0^+=t_0\le t \le t_1^-$ as follows.
We let evolve the  $N-M_\delta$ particles in $\und x^{\delta,-}(0^+)$ as in the basic process, namely as independent Brownian motions with non local branching.  If there is $\tau\in(0,t_1)$ such that $| \und x^{\delta,-}(\tau)|=N$, then $ \und x^{\delta,-}(t)$ for  $t>\tau$ is defined as independent Brownian motions without branching.  By definition
$N-M_\delta\le |\und x_\delta^-(t^-_1)|\le N$.  Define $ \und x^{\delta,-}(t_1^+)$ as the configuration obtained from $\und x_\delta^-(t^-_1)$ by deleting the $|\und x_\delta^-(t^-_1)|-(N-M_\delta)$ leftmost  particles. Thus $| \und x^{\delta,-}(t_1^+)|=N-M_\delta$ and so we can iterate the definition.

\end{itemize}

\noindent
It follows from their definition that $\und x^{\delta,\pm}(t)$ can be realized as subsets of the basic process $\und y(t)$.

\subsection{Stochastic inequalities}
\label{subsec3}

We will construct couplings to prove:

	\begin{theorem}
\label{thm1}
For each positive real number $\delta$
there is a coupling of the two processes $\und x(t)$ and $\und x^{\delta,+}(t)$ so that
	\begin{align}
	\label{2.3upper}
  \pi_{\und x(t)}\preccurlyeq  \pi_{\und x^{\delta,+}(t)},\qquad \text { for all }\,\,t
    \end{align}
There is also a a coupling of $\und x(t)$ and $\und x^{\delta,-}(t)$ so that
	\begin{align}
	\label{2.3lower}
 \pi_{\und x^{\delta,-}(t)}\preccurlyeq  \pi_{\und x(t)}\qquad \text { for all }\,\,t
\end{align}
\end{theorem}

We fix   $\delta>0$ and, to have lighter notations,  we will sometimes omit the dependence on $\delta$.  We
prove the lower bound in Subsection \ref{subse3.1} and  the upper bound in Subsection \ref{subse3.2}.

\subsection{Lower bound.}
\label{subse3.1}
We first introduce  an auxiliary process $\und z(t)$ of $N$ particles colored in red and blue in such a way that
the blue ones  have the same law as  $\und x^{\delta, -}(t)$.

The initial configuration  $\und z(0)$ is made by $N$ independent copies of variables with distribution $\rho_0(y)dy$; we paint in red the $M_\delta$ leftmost particles  and in blue the others.
To each blue particle we associate an independent exponential clock of parameter 1. When the clock rings (say at time $t$ for a blue particle at $x$) if there are no red particles we do nothing, otherwise we delete the rightmost red particle and we create  a new blue  particle at a position $x+Z$ where $Z$ here and in the sequel is a variable with  distribution $p(0,z)dz$.
In between branching times  the particles move as independent brownian motions.

At the times $k\delta$ we do a repainting:  let $m_k$ be the number of red particles at time $t_k^-=(k \delta)^-$, $k \ge 1$. By definition $0\le m_k\le M_\delta$. We then paint in red the $M_\delta-m_k$ leftmost blue particles so that at time $t_k^+=(k \delta)^+$ the number of red particles is again $M_\delta$.
Obviously the blue particles in the process $\und z(t)$ have the same law as  $\und x^{\delta, -}(t)$.

\medskip

\noindent {\bf {Coupling  the labelled true and auxiliary processes. }}
We are going to couple the labelled
processes $\und x(t)=\big(x_1(t),..,x_N(t)\big)$ and  $\und z(t)=\big(z_1(t),..,z_N(t)\big)$ calling $\und \si(t)=\big(\si_1(t),..,\si_N(t)\big)$, $\si_i(t)\in \{R,B\}$  the color of   $z_i(t)$.   Thus the coupled process is a process in $\mathbb R^N\times \mathbb R^N \times \{R,B\}^N$ whose elements are denoted by $(\und x,\und z,\und \si)$ and such that the marginal over $\und x$ is the true process $\und x(t)$ while the marginal law of $(\und z,\und \si)$ has the same law as the above auxiliary process $ \und z(t)$. The main point will be to  construct the coupling in such a way that at all times it is in the set
     \begin{equation}
      \label{4.1}
\mathcal X = \{(\und x,\und z,\und \si):
 x_i \ge z_i, 1\le i\le N\}
     \end{equation}
This will prove a stronger version of the lower bound \eqref{2.3lower} with the inequality holding for all $t$ (and not necessarily $t=k\delta$).  Indeed by the definition of the  process $\und z(t)$ the blue particles process is the process $\und x^{\delta,-}(t)$ and if $\si_i(t)=B$ then
$x_i^{\delta,-}(t) = z_i(t) \le x_i(t)$.  Thus the lower bound will be proved once we construct a coupling with values in $\mathcal X$ which we do next.

$\blacktriangleright$
At time 0 we set $x_i(0) = z_i(0)$, $i=1,..,N$, and define $\und \si(0)$ so that
 the $M_\delta$ leftmost particles of $\und z(0)$ are red while the others are blue.

$\blacktriangleright$ To each $i$ we associate independent exponential clocks of intensity 1 and  define the coupling inductively
between successive ``special times'' where a special time is either a time  when a  clock rings or  a time in $t_k$, $k\ge 0$.

$\blacktriangleright$  Let $t$ be a special time and suppose by induction that at time $t^+$ the coupled process is in $\mathcal X$.  We define the evolution till the successive special time so that the variables $x_i(s)$, $i=1,..,N$, are independent Brownian motions while each $z_i(s)$ has the same increments as $x_i(s)$ and keeps its color.  Hence in such a time interval the process is always in $\mathcal X$.

$\blacktriangleright$ Let the special time be $t_k$ and suppose by induction that at time $t_k^-$ the process is in $\mathcal X$.  We set $x_i(t_k^+)=x_i(t_k^-)$, $z_i(t_k^+)=z_i(t_k^-)$, $i=1,..,N$ and change only the colors in agreement with the definition of the $\und z$ process: namely we change from blue to red the color of the $M_\delta-m_k$ leftmost blue particles (recall that $m_k$ is the number of red particles at time $t_k^-$, $0\le m_k\le M_\delta$).  Since positions and labels are unchanged at time $t_k^+$ the process is still in $\mathcal X$.

 $\blacktriangleright$ Let $t$ be a time when the $i$-th clock rings and suppose by induction that at time $t^-$ the process is in $\mathcal X$.  In agreement with the definition
 of the $\und x$-process
  a new   $x$-particle is created at $x_i(t) + Z$
 while the leftmost $x$-particle is deleted, if this is the same as the new particle then $\und x(t)$ is unchanged, otherwise
  the new particle gets the same label as the particle which is deleted, namely
   \begin{equation*}
  x_h(t^+) =  x_h(t^-),\, h\ne j,\quad x_j(t^+)= \max\{ x_i(t^-)+Z, x_j(t^-)\}
  \end{equation*}
where  $j$ is the label of the leftmost $x$-particle,  i.e. $x_j(t^-) =\dis{ \min_{n}  x_n(t^-)}$. If
$\si_i(t)=R$ (in agreement with the rules of the $\und z$-process) we set $\und z(t^+)=\und z(t^-)$ so that $z_j(t^+) = z_j(t^-)\le x_j(t^-) \le x_j(t^+)$ and the coupled process stays in $\mathcal X$.
Suppose next that $\si_i(t)=B$ then a new blue particle in the $\und z$-process is created at position $z_i(t^-)+Z$, $Z$ the same as in the $x$-process, and simultaneously the rightmost red particle in $\und z(t^-)$ is deleted (if there was no red particle then $\und z(t^+)=\und z(t^-)$).  Call $k$ the label of the rightmost red particle in $\und z(t^-)$.  We set
   \begin{equation*}
  z_h(t^+) =  z_h(t^-),\, \si_h(t^+)=\si_h(t^-),\; h\notin \{ j,k\};\quad z_j(t^+)= z_i(t^-)+Z,\; \si_j(t^+)=B
  \end{equation*}
This completes the definition of $\und z(t^+)$ if $k=j$; if this is not the case, i.e.\ $j\ne k$, we complete the above definition by setting:
  \begin{equation*}
  z_k(t^+) =  z_j(t^-),\, \si_k(t^+)=\si_j(t^-)
  \end{equation*}
Since the coupled process is in $\mathcal X$, we have $z_j(t^+)= z_i(t^-)+Z\le  x_i(t^-)+Z
= x_i(t^+)$.  Moreover if  $j\ne k$, then
 \begin{equation*}
  z_k(t^+) =  z_j(t^-)\le x_j(t^-) \le x_k(t^-) = x_k(t^+)
  \end{equation*}

\subsection{Upper bound.}
\label{subse3.2}
We will define an auxiliary process $\und z(t)$ of red and blue colored particles in such a way that its marginal neglecting  the colors  has the law of $\und x^{\delta,+}(t)$.
 We will couple  $\und z(t)$ with the true process $\und x(t)$ in such a way that  $\und x(t) \preccurlyeq \und z^{(B)}(t)$, where $\und z^{(B)}(t)$ denotes the subset of $\und z(t)$ made of its blue particles;  this will prove the upper bound
 $\und x(t) \preccurlyeq\und x^{\delta,+}(t)$.

{\bf Definition.} The  auxiliary process $\und z(t)$.

\noindent
We define inductively $\und z(t)$ in the time intervals $[t_k^+,t_{k+1}^+]$.  We suppose by induction that   $|\und z(t_k^+)|=N$ and
 $\und z(t_k^+)=\und z^{(B)}(t_k^+)$, namely that at time $t^+_k$ there are $N$ particles and they are all blue, initially
$\und z(0)=\und z^{(B)}(0)=\und x(0)$.  We define $\und z^{(B)}(t), t\in [t_k^+,t_{k+1}^-]$ to have the law of the true process and at each time when a blue particle is deleted a red particle is created at the same place which then evolves as the basic process independently of all the other particles, the descendants of a red particle being all red.  Finally, to complete the induction step, we define $\und z(t_{k+1}^+)$ by retaining in $\und z(t_{k+1}^-)$ only the $N$ rightmost particles (independently of their color) and deleting all the others, we then paint in blue all those left.
Thus the process of red and blue particles together  has the law of  $\und x^{\delta,+}(t)$ while in each interval $[t_k^+,t_{k+1}^-]$,  $\und z^{(B)}(t)$ has the law of the true process.

Suppose now that $Q_{\und x,\und x'}$ is a coupling
 of the true processes    $\und x(t)$ and   $\und x'(t)$ starting from $\und x$ and $\und x'$, then we can construct a coupling of  $\und x(t)$ and $\und x^{\delta,+}(t)=\und z(t)$ by using in each time  interval $[t_k^+,t_{k+1}^-]$ the coupling $Q_{\und x(t^+_k),\und z(t^+_k)}$.
In
Lemma \ref{le2} below  we prove that $Q_{\und x,\und x'}$ can be chosen in such a way that  $Q_{\und x,\und x'}[\pi_{\und x(t)} \preccurlyeq \pi_{\und x'}(t)]=1$ for all $t> 0$ if $\pi_{\und x} \preccurlyeq \pi_{\und x'}$, from
this   the upper bound follows. The proof of Lemma \ref{le2}   uses the following elementary lemma.

	\begin{lemma}
	\label{le}
Let $\und x$ and $\und z$ $\in\mathbb R^N$, then the following statements are equivalent.
	\begin{description}
\item [S1.] $\pi_{\und x} \preccurlyeq \pi_{\und z}$
\item [S2.] There are $i_1,..,i_N\in\{1,..,N\}$ and  $j_1,..,j_N\in\{1,..,N\}$
such that $x_{i_\ell}\le z_{j_\ell}$, for all $\ell=1,..,N$
\item [S3.]  Call $i'_1$ the label  of the leftmost $x$-particle, $i'_2$ of the second one ...,$i'_N$ of the last one  so that
  $x_{i'_1}\le x'_{i_2} \le\cdots \le x_{i'_N}$.  Do the same for   $\und z$ calling $j'_1,..,j'_N$ the corresponding labels, so that
  $z_{j'_1}\le z_{j'_2} \le\cdots \le z_{j'_N}$. Then $x_{i'_\ell}\le z_{j'_\ell}$, for all $\ell=1,..,N$.
	\end{description}
	\end{lemma}

	\begin{lemma}
	\label{le2}
Let $\und x$ and $\und x'$ be in $\mathbb R^N$ and let
$\pi_{\und x} \preccurlyeq \pi_{\und x'}$.  Then there is a coupling $P$
of the true processes starting from $\und x$ and $\und x'$ such that $P[\pi_{\und x(t)} \preccurlyeq \pi_{\und x'}(t)]=1$ for all $t\ge 0$.

	\end{lemma}

\noindent
{\bf Proof.}  We introduce $N$ clocks, $i=1,..,N$, which ring  independently at rate 1. The idea of the coupling is that when a clock ring, say clock $i$, then the $i$-th particle counting from the left in the $\und x$ and the $i$-th particle counting from the left  in the $\und x'$ process branch in the same way.  More precisely we define iteratively the process in the time intervals $[s^+_k,s^+_{k+1}]$ when the clocks ring.  We suppose inductively that
$\pi_{\und x(s^+_k)} \preccurlyeq \pi_{\und x'}(s^+_k)$, so that by Lemma \ref{le} if we label the particles counting from the left (i.e.\ $x_i\le x_{i+1}, i=1,..,N-1$), then $x_i(s^+_k)\le  x'_i(s^+_k)$ for all $i=1,..,N$.
We define the coupled process in the time interval $[s^+_k,s^-_{k+1}]$ by letting $x_i(s)$, $i=1,..,N$, be independent Brownian motions and $x'_i(s)$ have the same increments as $x_i(s)$.  Then evidently $x_i(s^-_{k+1})\le  x'_i(s^-_{k+1})$ for all $i=1,..,N$, so that by Lemma \ref{le} $\pi_{\und x(s^-_{k+1})} \preccurlyeq \pi_{\und x'}(s^-_{k+1})$. We can relabel the particles  counting from the left, as we did before so that $x_j(s^-_{k+1})\le  x'_j(s^-_{k+1})$ for all $j=1,..,N$.  Let us now suppose that it is the $i$-th clock which rings.  We then set, calling for brevity $t^{\pm}=s^{\pm}_{k+1}$: $x_j(t^+)=x_j(t^-)$, $x'_j(t^+)=x'_j(t^-)$, $j=2,..,N$ and
     $$
     x_1(t^+) = \max\{ x_i(t^-) + Z, x_1(t^-)\},\quad x'_1(t^+) =
      \max\{ x'_i(t^-) + Z, x'_1(t^-)\}
      $$
 Then  $x_j(t^+)\le x'_j(t^+)$, $j=1,..,N$ and the induction hypothesis is verified.  \qed

\section{Proof of Theorem \ref{thmDD2.1}}
\label{sec4}

In this section we prove Theorem \ref{thmDD2.1} using  Theorem \ref{thmDD2.3} and
Theorem \ref{thmDD44.1} below  (proved later in Section \ref{sec:DDF.7}).
Theorem \ref{thmDD44.1}  states that
the stochastic upper and lower barriers are
with ``large'' probability ``close'' to the corresponding deterministic barriers  for large $N$.
Closeness is quantified using the following semi-norms:

\noindent {\bf Semi-norms.} {\em Let $\mu$ and $\nu$ be positive, finite
measures on $\mathbb R$ and $\mathcal I_N$ the partition of $\mathbb R$ into intervals $I=[kN^{-\beta},(k+1)N^{-\beta})$, $k\in \mathbb Z$, (we will eventually fix $ \beta = \frac 1{12}$).  We then define } (below $\mathcal A \subset \mathcal I_N$)
	 \begin{equation}
	 \label{DD44.1}
\|\mu-\nu\|_{\mathcal I_N} := \sum_{I \in \mathcal I_N } |\mu(I)-\nu(I)|,\quad
\|\mu-\nu\|_{\mathcal A}=\sum_{I \in \mathcal A} |\mu(I)-\nu(I)|
   \end{equation}
.

The semi-norm $\|\mu-\nu\|_{\mathcal I_N}$ is the variational distance of the marginal distributions of $\mu$ and $\nu$ on the $\si$-algebra generated
by $\mathcal I_N$. It is also the $L^1$-norm of coarse grained versions of $\mu$ and $\nu$ on the scale $N^{-\beta}$ which are defined as
	 \begin{equation}
	 \label{DD44.2}
\phi(r) = \sum_{I\in \mathcal I_N}\frac{ \nu[I] }{|I|}\mathbf 1_{r\in I},\quad \psi(r) = \sum_{I\in \mathcal I_N}\frac{ \mu[I] }{|I|}\mathbf 1_{r\in I}
   \end{equation}
Indeed
\eqref{DD44.1}	 can be obviously written as
  \begin{equation}
	 \label{DD44.3}
\|\nu -\mu \|_{\mathcal I_N} = \int dr\, |\phi(r)-\psi(r)|
   \end{equation}

The semi-norms control mass transport order (see \eqref{DD2.2}):

\begin{lemma}

With the above notation
	\begin{eqnarray}
	\label{DD2.4n}
\sup_{r\in \mathbb R}\Big |\mu\Big[[r,\infty)\Big] - \nu\Big[[r,\infty)\Big] \Big| &\le& \|\nu -\mu \|_{\mathcal I_N} +
\sup_{I \in \mathcal I_N} \{ \mu[I]+\nu[I]\} \\ &\le& 2
\|\nu -\mu \|_{\mathcal I_N} +\sup_{I \in \mathcal I_N}\nu[I] \nn
     \end{eqnarray}

\end{lemma}

\noindent
{\bf Proof.}  Fix $r\in \mathbb R$ and call $I_r$ the interval in $\mathcal I_N$ which contains $r$; write $I>r$ for the intervals to the right of $I_r$.  Then
	\begin{eqnarray*}
\Big|\mu \big[[r,\infty)\big] - \nu \big[[r,\infty)\big]\Big| &\le& \sum_{I>r}|\nu[I] -\mu[I] |  +| \int _{I_r}
\mu(dr')
\mathbf 1_{r'\ge r} - \int _{I_r}
\nu(dr') \mathbf 1_{r'\ge r}| \\ &\le&  \|\nu -\mu \|_{\mathcal I_N} + \mu[I_r]+\nu[I_r]
     \end{eqnarray*}
hence the first inequality in \eqref{DD2.4n}.  The second one follows because
$ \mu[I] \le \nu[I] +    \|\nu -\mu \|_{\mathcal I_N}$.  \qed

     \smallskip
We will use \eqref{DD2.4n} with $\nu$ a measure with bounded density with respect to Lebesgue so that $\nu(I)\le c|I| =cN^{-\beta}$.
Therefore the last term in \eqref{DD2.4n} will be negligible.
In fact
we will use the semi-norms to compare the counting measures
$\frac 1N\pi_{\und x^{\delta,\pm}(t)}(dr)$ associated to the stochastic barriers and the measures
$\rho^{\delta,\pm}(r,t)dr$  associated to the deterministic barriers, observing that $\rho^{\delta \pm}(r,t)$ are uniformly bounded. As explained in Section \ref{sec:DD2} we will work in the following context.

We fix arbitrarily $T>0$,  $ n\in \mathbb N$, call $\delta = 2^{-n}T$, $t_k= k \delta$, $k=0,..,2^{n}$,
 and shorthand
	\begin{align}
	\label{DD2.4n.1}
\und x^{\delta,\pm}_k = \und x^{\delta,\pm}(t_k^+),\quad
\rho^{\delta,\pm}_k(dr) = \rho^{\delta,\pm}(r,t_k^+)dr; \quad
 k=0,..,2^n
     \end{align}

We will prove in Section \ref{sec:DDF.7} the following Theorem.

\begin{theorem}
\label{thmDD44.1}
Let $\beta=\frac 1{12}$,  $\alpha^*= \frac 1{12}$ and
\begin{equation}
\label {DD44.5bis} \zeta_{N,T} = e^T\delta^{-1/2} N^{-\alpha^*}+e^{2T}\frac {T^{N^{\alpha^*/2}}}{N^{\alpha^*/2}!}
\end{equation}
 Then  there
is $c^*$ so that the set
 	 \begin{equation}
	 \label{DD44.5}
\|\frac 1N\pi_{\und x^{\delta,\pm}_k}-\rho^{\delta,\pm}_k\|_{\mathcal I_N}\le c^*\zeta_{N,T},\quad \text{for all $k\le K=2^n$}
   \end{equation}
has $P^{(N)}$-probability  $\ge 1 -c^* 2^{n}  N^{-\alpha^*}$.

\end{theorem}

\noindent {\bf Proof of Theorem \ref{thmDD2.1}} (pending the validity of Theorem \ref{thmDD44.1}).

By Theorem \ref{thmDD44.1} and  \eqref{DD2.4n}, by \eqref{DDD3.3.0} and \eqref{DD3.24}
$\rho^{\delta,\pm}(r,t) \le c'\equiv e^T$,
	\begin{eqnarray*}
P^{(N)}\Big[\bigcap_{k}\big \{\sup_{r\in \mathbb R}\Big |\frac 1N\pi_{\und x^{\delta,\pm}_k}\Big[[r,\infty)\Big] - \rho^{\delta,\pm}_k\Big[[r,\infty)\Big] \Big| &\le&   2 c^*\zeta_{N,T}
  +c' N^{-\beta}\big \}\Big] \\
  \ge 1 -c^* 2^{n}  N^{-\alpha^*}
  \nn
     \end{eqnarray*}
With reference to the upper barrier we thus have
	 \begin{equation*}
P^{(N)}\Big[\bigcap_{k}\{\frac 1N\pi_{\und x^{\delta,+}_k} \preccurlyeq \rho^{\delta,+}_k,\;\;
  \text{modulo $\eps$}\}\Big] \ge 1 -c^* 2^{n}  N^{-\alpha^*},\;
   \eps:=2 c^*\zeta_{N,T}
  +c' N^{-\beta}
\end{equation*}
  Call $\mathcal P$ the law of the coupling defined in
 Theorem \ref{thm1}, then
	\begin{eqnarray*}
P^{(N)}\Big[\bigcap_{k}\{\frac 1N\pi_{\und x^{\delta,+}_k} \preccurlyeq \rho^{\delta,+}_k,\;\;
  \text{modulo $\eps$}\}\Big] = \mathcal P\Big[\bigcap_{k}\{\frac 1N\pi_{\und x^{\delta,+}_k} \preccurlyeq \rho^{\delta,+}_k,\;\;
  \text{modulo $\eps$}\}\Big]\\ \ge
  \mathcal P\Big[\bigcap_{k}\{\frac 1N\pi_{\und x_k} \preccurlyeq \rho^{\delta,+}_k,\;\;
  \text{modulo $\eps$}\}\Big] = P^{(N)}\Big[\bigcap_{k}\{\frac 1N\pi_{\und x_k} \preccurlyeq \rho^{\delta,+}_k,\;\;
  \text{modulo $\eps$}\}\Big]
     \end{eqnarray*}
hence
	\begin{eqnarray*}
P^{(N)}\Big[\bigcap_{k}\{\frac 1N\pi_{\und x_k} \preccurlyeq \rho^{\delta,+}_k,\;\;
  \text{modulo $\eps$}\}\Big ]\ge 1-c^* 2^n N^{-\alpha^*}
     \end{eqnarray*}
An analogous argument holds for the lower barriers so that
	\begin{eqnarray}
	\label{DF.1}
P^{(N)}\Big[\bigcap_{k}\{\rho^{\delta,-}_k\preccurlyeq
\frac 1N\pi_{\und x_k} \preccurlyeq \rho^{\delta,+}_k,\;\;
  \text{modulo $\eps$}\}\Big ]\ge 1-2c^* 2^{n} N^{-\alpha^*}
     \end{eqnarray}
On the other hand by \eqref{DD2.9b}
  	 \begin{equation}
	 \label{DD2.9?}
\rho^{\delta,+}(r,t)dr \preccurlyeq \rho^{\delta,-}(r,t)dr\quad \text{modulo $c \delta$}
   \end{equation}
hence by \eqref{DD2.9}
 	 \begin{equation}
	 \label{DD2.9??}
\rho^{\delta,+}(r,t)dr \preccurlyeq u(r,t)dr\quad \text{modulo $c \delta$}
   \end{equation}
Since an analogous argument holds for the lower barriers we get from \eqref{DF.1}
	\begin{eqnarray}
	\label{DF.1.1}
P^{(N)}\Big[u(r,t_k)dr\preccurlyeq
\frac 1N\pi_{\und x_k} \preccurlyeq u(r,t_k)dr,\;\;
  \text{modulo $\eps+c\delta$}\Big ]\ge 1-2c 2^{n} N^{-\alpha^*}
     \end{eqnarray}
  We then get  Theorem \ref{thmDD2.1} recalling that  $\eps:=2 c^*\zeta_{N,T}
  +c' N^{-\beta} $, with $\zeta_{N,T}$ defined in \eqref{DD44.5bis}.

\section{Properties of the basic process}
\label{sec5}

 The proof of Theorem \ref{thmDD44.1} splits into two parts, the analysis of $\und x^{\delta,\pm}(t)$ and $N\rho^{\delta,\pm}(\cdot,t)$ for
$t\in (k\delta,(k+1)\delta)$ and the analysis of the transition from time
$((k+1)\delta)^-$ to $((k+1)\delta)^+$.  The former is common to lower and upper barriers, in the second we have to distinguish between upper and lower barriers.
In the time interval $ (k\delta,(k+1)\delta)$ the process is without deaths, it is therefore the basic process defined in Section \ref{sec:DD2}. In the next theorem we will prove that in the average the basic process behaves as the deterministic free evolution of Subsection \ref{subseDD3.1}.  In an appendix we will use this to prove some estimates which are used in Section \ref{sec:DDF.7} to prove that the stochastic barriers converge to the deterministic barrier in the continuum limit $N\to \infty$.

\smallskip

Recall that $P^{\und y_0}$ denotes the law of
the basic process starting from the configuration $\und y_0$.
Calling $P^x_t(d\und y)$ the law of $\und y(t)$ starting from the configuration $\und y_0$ consisting of a single particle at $x$,
 we define the averaged counting  measure at time $t$ starting from $x$ as
	 \begin{equation}
	 \label{DD5.1}
\la^x_t(dr):= \int P^x_t(d\und y)\pi_{\und y}(dr)
   \end{equation}
The next theorem is a key step in the comparison between stochastic and deterministic barriers.
Recall from Section  \ref{sec:DD3} that we have called $L$ the adjoint of $L^*$ defined in \eqref{DD3.5}.

\begin{theorem}
\label{thmDD5.1}
With the above notation  $\la^x_t(dr)= e^{Lt}(x,r)dr$.

\end{theorem}

\noindent
{\bf Proof.}
Let $\phi$ be a smooth test function and call
	 \begin{equation}
	 \label{DD5.7}
h(x,t)= \int \la^x_t(dr)\phi(r)
   \end{equation}
By translation invariance for any $a\in \mathbb R$
	 \begin{equation*}
	\nn
h(x+a ,t)= \int  \la^x_t(dr)\phi(r-a)
   \end{equation*}
which proves that $h(x,t)$ is a smooth function of $x$.  It then follows that
	 \begin{equation*}
 \nn
\frac{\partial}{\partial t} h(x,t)= \frac 12 \frac{\partial^2}{\partial x^2}h(x,t)
+ \int dy \; p(x,y) h(y,t)
 \end{equation*}
with $h(x,0)= \phi(x)$.  Thus $h(x,t)= e^{L t}\phi(x)$.
By \eqref{DD5.7} $ \int  \la^x_t(dr)\phi(r) =  e^{L t}\phi(x)$
hence $\la^x_t(dr)=e^{Lt}(x,r)dr$.  \qed
\smallskip

Observe that since the branchings are independent:
        	 \begin{equation}
	\nn
 \int P^{\und y_0}_t(d\und y)\phi(\und y) =\int\{ \prod_{x\in \und y_0}P^{x }_t(d\und y^{(x)})\}\phi\Big(\bigcup_{x\in \und y(0)}\und y^{(x)}\Big)
   \end{equation}
and in particular
    \begin{equation}
	 \label{DD5.22}
 \int P^{\und y_0}_t(d\und y) \int \pi_{\und y} (dr)f(r) = \sum_{x\in \und y_0}\int P^{x }_t(d\und y)\int \pi_{\und y}(dr ) f (r)
   \end{equation}

\section{Continuum limit of the stochastic barriers }
\label{sec:DDF.7}

To prove Theorem \ref{thmDD44.1} (to which we refer for notation) we must find a ``good'' set $\mathcal X^{\rm good}$ of large probability where the semi-norms
$\| \frac 1N
\pi_{\und x^{\delta,\pm}_{k+1}}
-\rho^{\delta,\pm}_{k+1}\|_{\mathcal I_N}$ are small.

\subsection {The good set}

$\mathcal X^{\rm good}$ is the intersection of a good set for the upper barrier and a good set for the lower barrier, which are both intersections of four good sets, thus
	 \begin{equation}
	 \label{DDF.1.0}
\mathcal X^{\rm good} = \Big\{ \bigcap_{i=1}^4\mathcal X^{+}_i \Big\} \cap
\Big\{\bigcap_{i=1}^4\mathcal X^{-}_i \Big\}
   \end{equation}
All sets are tacitly defined on the same space which is the space where  the basic process $\und y(t), t\ge 0$ is realized.  They will be defined using parameters which  should satisfy the conditions:
   	 \begin{equation}
	 \label{DD5.30}
 \frac{1+(\beta+b)}2 < \alpha_1 < 1-(\beta+b),\;\;\beta+b \in (0, \frac 1{6}], \quad
 \alpha_0 > \frac 12
   \end{equation}
Our specific choice is
	 \begin{equation}
	 \label{DDF.0}
\beta=b=\frac1{12},\; \alpha_1=\alpha_0= \frac 23,\;
   \end{equation}

We are now ready to define the good set:

$\bullet$ The first good set is:
 	 \begin{equation}
	 \label{DDF.1}
\mathcal X^{\pm}_1:= \bigcap_{k=1}^K  \{\und x_k^{\pm} \cap [-N^b,N^b]^c = \emptyset\}, \quad K=2^n
   \end{equation}
To define the second good set we need the following notation:
	 \begin{equation}
	 \label{DDF.3}
\mathcal I_N= \mathcal I'_N\cup \mathcal I''_N,\quad
\mathcal I'_N=\{I\in \mathcal I_N: I \cap [-N^b,N^b]\ne \emptyset\}
   \end{equation}
$\mathcal I''_N$ being the set of all $I \in \mathcal I_N$ which have empty intersection with   $[-N^b,N^b]$.

$\bullet$  Calling $t_k^-=(k\delta)^-$, the second good set is then:
 	 \begin{eqnarray}
	 \label{DDF.5}
&&\mathcal X^{\pm}_2 := \bigcap_{k=1}^K \Big\{ \| \frac 1N \pi_{\und x^{\delta,\pm}(t_k^-)}(dr)-\rho^{\delta,\pm}(r,t_k^-)(dr)\|_{\mathcal I'_N} \le   e^{\delta}
\| \frac 1N \pi_{\und x^{\delta,\pm}_{k-1}}-\rho^{\delta,\pm}_{k-1}\|_{\mathcal I_N}
\nn\\ &&  \hskip4cm+  c' \delta^{-1/2} N^{-\beta} + c''N^{\beta+b+\alpha_1-1}\Big\}
   \end{eqnarray}

$\bullet$ The third good set $\mathcal X^{\pm}_3$ involves the values $n^{\pm}_k$ of the number of particles in the stochastic barriers at the times $t_k^-:= (k\delta)^-$, namely $ n^{\pm}_k := \pi_{\und x^{\delta,\pm}(t_k^-)}[\mathbb R]$:
 	 \begin{equation}
	 \label{DDF.6}
\mathcal X^+_3:= \bigcap_{k=1}^K \{|n^+_k -e^{\delta}N| \le N^{\alpha_1}\},\quad
\mathcal X^-_3:= \bigcap_{k=1}^K
 \{|n^-_k -[N-e^{\delta}N^{\alpha_0}]| < N^{\alpha_0}\}
   \end{equation}
   
$\bullet$  The fourth good set $\mathcal X^{\pm}_4$ involves what happens at time 0:
 \begin{equation}
	 \label{DDF.7}
\mathcal X^{\pm}_4:=\{
 \| \frac 1N \pi_{\und x^{\delta,\pm}_0}-\rho^{\delta,\pm}(0,r)dr\|_{\mathcal I_N} \le
 c'''N^{\beta+\alpha_1-1}\}
    \end{equation}

By choosing the parameters
$c'\dots c'''$  sufficiently large we will prove that the probability of the good sets will go to 1 as $N\to \infty$.  Thus the proof of Theorem \ref{thmDD44.1} splits into two parts, in the first one we show that in the good set the semi-norms are close and a second part where we prove that the probability of not being in the good set vanishes as $N\to \infty$.

\subsection{Bound in the good set}

We are going to prove that  with the choice \eqref{DDF.0} of the parameters in the good set
$\mathcal X^{\rm good}$:
 	 \begin{eqnarray}
	 \label{DDF.17}
&&\| \frac 1N \pi_{\und x^{\delta,+}_k}-\rho^{\delta,+}_k(dr)\|_{\mathcal I_N} \le c\Big(
 e^{T} \delta^{-1/2} N^{-1/12} + e^{2T}\frac {T^{N^{1/24}}}{N^{1/24}!}\Big)
      \end{eqnarray}
$c$ a suitable constant.  The bound  for the lower barriers is the same, its proof is similar and omitted.
We fix in the sequel $k\in \{1,..,2^n\}$, the bounds will be uniform in $k$.

Our first task is:
$$\text{ bound }\|  \pi_{\und x^{\delta,+}_k}-N\rho^{\delta,+}_k(dr)\|_{\mathcal I_N} \text{ in terms of }
\|  \pi_{\und x^{\delta,+}(t_k^-)}-N\rho^{\delta,+}(r,t_k^-)(dr)\|_{\mathcal I_N}$$
To have lighter notation we write $d\mu''= \pi_{\und x^{\delta,+}_k}(dr)$, $d\nu''=N\rho^{\delta,+}_k(dr)$ and
$d\mu= \pi_{\und x^{\delta,+}(t_k^-)}(dr)$, $d\nu=N\rho^{\delta,+}(r,t_k^-)dr$.
$\mu''$ is obtained from $\mu$ by cutting on the left a mass $e^{\delta}N - N + \theta$, by \eqref{DDF.6} $|\theta| \le N^{\alpha_1}$.   Instead $\nu''$ is obtained from $\nu$
by cutting on the left a mass $\Theta=e^{\delta}N - N$.

If $\theta=0$ we are cutting the same mass from $\mu$ and $\nu$; suppose $\mu$ has density $f$, $\nu$ has density $g$ and call $f''$ and $g''$ the densities of $\mu''$ and $\nu''$.
Since we are cutting mass from the left the $L^1$-norm of $f''-g''$ is not larger than that of $f-g$, see for instance  Proposition 5.2 in \cite{CDGP}.  The whole point will be to prove that the same property holds for the semi-norms, see Proposition \ref{lemmaDDF7.1} below.  In general however $\theta$ will not be zero, suppose for the sake of definiteness that $\theta \ge 0$.  The following is to reduce to the case $\theta=0$.
Let then $\la''$ be obtained from $\mu$
by cutting on the left a mass $\Theta=e^{\delta}N - N$, so that $\mu''=\la'' + \rho$ where $\rho$ is a positive measure with mass $\theta$.  Then
   $$
 \| \mu''-\la''\|_{\mathcal I_N} = \int \rho(dr) = \theta
 $$
and therefore
	 \begin{equation}
	 \label{DDF.8}
 \| \mu''-\nu''\|_{\mathcal I_N} \le  \| \mu''-\la''\|_{\mathcal I_N} +  \| \la''-\nu''\|_{\mathcal I_N} \le \theta +  \| \la''-\nu''\|_{\mathcal I_N}
    \end{equation}

    \medskip
    The next proposition improves an analogous statement in Lemma 11.6 of \cite{CDGP}:

\begin{proposition}
\label{lemmaDDF7.1}
Let $\la''$ and $\nu''$ be obtained from the positive measures  $\mu$ and $\nu$
by cutting on the left a mass $\Theta$ smaller than $\mu[\mathbb R]$ and $\nu[\mathbb R]$, then
	 \begin{equation}
	 \label{DDF.9}
 \| \la''-\nu''\|_{\mathcal I_N} \le \| \mu-\nu\|_{\mathcal I_N}
    \end{equation}

\end{proposition}

\noindent
{\bf Proof.}
     Let $\psi''$ and $\psi$ be
the coarse grained versions of $d\la''$ and $d\mu$; $\phi''$ and $\phi$
those relative to $d\nu''$ and $d\nu$, see \eqref{DD44.2} for notation.
Let $R_\nu$ and  $R_\mu$ be such that
   	 \begin{equation}
	 \label{DDF.10}
\int_{-\infty}^{R_\nu}dr\, \phi(r) = \Theta,\quad
\int_{-\infty}^{R_\mu}dr\, \psi(r) = \Theta
   \end{equation}
Observe that it is not necessarily true that $\mu[(-\infty,R_\mu]] = \Theta$
and $\nu[(-\infty,R_\mu]] = \Theta$.

Suppose for the sake of definiteness that $R_\mu \ge R_\nu$.  Call
   	 \begin{equation}
	 \label{DDF.11}
\phi'(r) = \phi(r) \mathbf 1_{r\ge R_\nu},\quad
\psi'(r) = \psi(r) \mathbf 1_{r\ge R_\mu}
   \end{equation}
For    any $I \in \mathcal I_N$:
      	 \begin{equation}
	 \label{DDF.12}
 \int_I dr\, \phi'(r) =  \int_I dr\, \phi''(r),\quad
 \int_I dr\, \psi'(r) =  \int_I dr\, \psi''(r)
   \end{equation}
 Since the cutting operation does not increase the $L^1$ norm, see Proposition 5.2 in \cite{CDGP},
         	 \begin{equation}
	 \label{DDF.13}
\int dr\, |\phi'(r) -\psi'(r)| \le \int dr\, |\phi(r) -\psi(r)|= \| \mu-\nu\|_{\mathcal I_N}
   \end{equation}
 We have
    $$
    \int dr\, |\phi''(r) -\psi''(r)| \le \int dr\, |\phi'(r) -\psi'(r)|
    $$
    so that, by \eqref{DDF.12}  and  \eqref{DDF.12}:
         	 \begin{equation*}
  \|\nu'' -\la'\|_{\mathcal I_N}=\int dr\, |\phi''(r) -\psi''(r)| \le \int dr\, |\phi'(r) -\psi'(r)| \le\|\nu-\mu\|_{\mathcal I_N}
    \end{equation*}
\qed

By \eqref{DDF.8} and Proposition \ref{lemmaDDF7.1} in $\mathcal X_3^{+}$
   \begin{equation}
	 \label{DDF.144}
  \| \frac 1N \pi_{\und x^{\delta,+}_k}-\rho^{\delta,+}_k\|_{\mathcal I_N} \le
     \| \frac 1N \pi_{\und x^{\delta,+}(t_k^-)}(dr)-\rho^{\delta,+}(r,t_k^-)dr
   \|_{\mathcal I_N} + N^{\alpha_1-1}
   \end{equation}

\begin{lemma}
\label{lemmaDDF7.2}
In $\mathcal X_1^+$
	 \begin{eqnarray}
	 \label{DDF.14}
 \| \frac 1N \pi_{\und x^{\delta,+}(t_k^-)}(dr)-\rho^{\delta,+}(r,t_k^-)dr\|_{\mathcal I_N} &\le&
    \| \frac 1N \pi_{\und x^{\delta,+}(t_k^-)}(dr)-\rho^{\delta,+}(r,t_k^-)dr\|_{\mathcal I'_N}\nn
    \\&+& c_b e^T \frac{T^{N^{b/2}}}{N^{b/2} !}
   \end{eqnarray}

\end{lemma}

\noindent
{\bf Proof.}
In $\mathcal X_1^+$ $\pi_{\und x^{\delta,+}_k}[I]=0$ for all  $I \in \mathcal I''_N$. Then, since $k\delta \le T$,
\eqref{DDF.14}  follows using \eqref{DD5.12}.  \qed

Then by \eqref{DDF.144}, \eqref{DDF.14} and \eqref{DDF.5}
	 \begin{eqnarray}
	 \label{DDF.15}
&&\| \frac 1N \pi_{\und x^{\delta,+}_k}-\rho^{\delta,+}_k\|_{\mathcal I_N} \le
e^{\delta}
\| \frac 1N \pi_{\und x^{\delta,+}_{k-1}}-\rho^{\delta,+}_{k-1}\|_{\mathcal I_N} + \eps_N \nn\\&&
\eps_N:= N^{\alpha_1-1} + c_b e^T \frac{T^{N^{b/2}}}{N^{b/2} !}
+ c' \delta^{-1/2} N^{-\beta} + c''N^{\beta+b+\alpha_1-1}\}
   \end{eqnarray}
   Hence by  \eqref{DDF.7} 	
	 \begin{eqnarray}
	 \label{DDF.16}
&&\| \frac 1N \pi_{\und x^{\delta,+}_k}-\rho^{\delta,+}_k\|_{\mathcal I_N} \le
   e^{\delta k} c'''N^{\beta+\alpha_1-1} + \sum_{h=0}^{k-1} e^{\delta h} \eps_N
      \end{eqnarray}
\eqref{DDF.17} is therefore proved with the choice \eqref{DDF.0} of the parameters.

\subsection{Probability of the good set}

The hardest part is to estimate the probability of the good sets $\mathcal X^{\pm}_2$.
In the case of  $\mathcal X^{+}_2$ we will exploit the fact that in each interval
$(t_{k-1},t_k)$ the upper barrier $x^{\delta,+}(t)$ is in law the same as the basic process $\und y(t)$.  This is however no longer true for  the lower barrier $x^{\delta,-}(t)$.  In fact to prove the stochastic inequality for the lower barrier we needed to stop the branching as soon as $|x^{\delta,-}(t)|=N$.  To deal with that we use the following ``trick''.
We first define a new lower stochastic barrier denoted by $\und z^{\delta,-}(t)$ which is defined like $\und x^{\delta,-}(t)$ but without stopping the branching when the number of particles becomes equal to $N$. We then define
	 \begin{eqnarray}
	 \label{DDF.5bis}
&&\mathcal Z^{-}_2 := \bigcap_{k=1}^K \Big\{ \| \frac 1N \pi_{\und z^{\delta,-}(t_k^-)}(dr)-\rho^{\delta,-}(r,t_k^-)(dr)\|_{\mathcal I'_N} \le   e^{\delta}
\| \frac 1N \pi_{\und z^{\delta,-}_{k-1}}-\rho^{\delta,-}_{k-1}\|_{\mathcal I_N}
\nn\\ &&  \hskip4cm+  c' \delta^{-1/2} N^{-\beta} + c''N^{\beta+b+\alpha_1-1}\Big\}
   \end{eqnarray}

\begin{lemma}
Recalling the definition \eqref{DDF.6} of $\mathcal X^-_3$,
	 \begin{equation}
	 \label{DDF.6.0}
\mathcal X^-_2 \cap \mathcal X^{-}_3 = \mathcal Z^-_2 \cap \mathcal X^{-}_3
   \end{equation}

\end{lemma}

\noindent {\bf Proof.}
In $\mathcal X^-_3$
 the total number of particles
in  $x^{\delta,-}(t_k^-)$ denoted by $n^-_k$ is  bounded by
   $$
n^-_k - [N-e^{\delta}N^{\alpha_0}] + [N-e^{\delta}N^{\alpha_0}] \le
 N^{\alpha_0} -e^{\delta}N^{\alpha_0} +N <N
   $$
 Then  in   $\mathcal X^-_3$ $x^{\delta,-}(t)= z^{\delta,-}(t)$ for  all $t$.  \qed

\smallskip
Thus the probability of $\mathcal X^{\rm good}$ is the same as that with  $\mathcal X^-_2$ replaced by $\mathcal Z^-_2$: in the sequel we will estimate the latter.
We obviously have:
	 \begin{eqnarray}
	 \label{DDF.56}
P^{(N)}[\mathcal X^{\rm good}]
&\ge& 1 - \Big(\{P^{(N)}[
(\mathcal X^{+}_1)^c] + P^{(N)}[
(\mathcal X^{-}_1)^c] \}\nn\\&+& \sum_{k=1}^{2^n}\{P^{(N)}[
(\mathcal X^{+}_2(k))^c] +P^{(N)}[
(\mathcal Z^{-}_2(k))^c]\} \nn
 \\&+& \sum_{k=1}^{2^n}\{P^{(N)}[
|n^+_k -e^{\delta}N| > N^{\alpha_1}] +P^{(N)}[
|n^-_k -[N-e^{\delta}N^{\alpha_0}]| > N^{\alpha_0}]\} \nn
 \\&+& P^{(N)}[\|
  \frac 1N \pi_{\und x_0}-\rho_0(r)dr\|_{\mathcal I_N} >
 c'''N^{\beta+\alpha_1-1}]  \nn
 \\&+&
 P^{(N)}[\|
  \frac 1N \pi_{\und x^{\delta,-}_0}-\rho^{\delta,-}(0,r)dr\|_{\mathcal I_N} >
 c'''N^{\beta+\alpha_1-1}]
  \Big)
   \end{eqnarray}
where $\mathcal X^{+}_2(k)$  are the sets in \eqref{DDF.5} involving $t_k^-$,  $\mathcal Z^{-}_2(k)$ is defined analogously.

\smallskip

By \eqref{DD5.13.1}
	 \begin{equation}
    \label{{DDF.17}}
P\Big[ \mathcal X^{\pm}_1\Big] \ge 1 - c_{b,T} N^{-1}
 \end{equation}
The second term on the right hand side of \eqref{DDF.56} will be bounded at the end of the subsection.
The third term in \eqref{DDF.56} is bounded by:
  	 \begin{equation}
	 \label{DDF.18}
 \le  2^n [c \delta N^{1-2\alpha_1}+c \delta N^{1-2\alpha_0}]=   2^n 2 c \delta N^{1-2\alpha_0}
   \end{equation}
having taken $\alpha_1=\alpha_0$. To prove \eqref{DDF.18} we use \eqref{DD5.25} with $\und y_0= \und x^{\delta,\pm}_{k-1}$.  For the upper barrier  $\la^{\und y_0}_\delta(\mathbb R) = N$ because  $\pi_{\und x^{\delta,+}_{k-1}}[\mathbb R]= e^{-\delta }N$. For the lower barrier
 $N - e^\delta N^{\alpha_0}$ because
  $\pi_{\und x^{\delta,-}_{k-1}}[\mathbb R]= e^{-\delta }N-N^{\alpha_0}$.

For the fourth term we use the Chebishev inequality with the squares to get
\begin{equation*}
 P^{(N)}[|
   \pi_{\und x_0}[I]-N\int_I dr\,\rho_0(r)|  \le
  N^{ \alpha_1}] \ge 1- c N^{1-2\alpha_1} N^{-\beta}
   \end{equation*}
Recall in fact that $\und x_0$ is the configuration obtained by taking $N$ independent copies
distributed as $\rho_0(r)dr$. Then
\begin{equation}
	 \label{DDF.19}
 P^{(N)}\Big[\bigcap _I\;\{ |
   \pi_{\und x_0}[I]-N\int_I dr\,\rho_0(r)|  \le
  N^{ \alpha_1}\}\Big] \ge 1- c' N^{\beta}c N^{1-2\alpha_1} N^{-\beta}
   \end{equation}
 because $\rho_0$ has compact support so that the number of intervals $I$ to consider is   $\le c'N^{\beta}$.  In the   set on the left hand side of \eqref{DDF.19}
    $$
     |
  \sum_I |\frac 1N \pi_{\und x_0}[I]-\int_I dr\,\rho_0(r)| \le c' N^{\beta} N^{ \alpha_1-1}
  $$
  so that
  \begin{equation}
	 \label{DDF.19.0}
  P^{(N)}\Big[\|
  \frac 1N \pi_{\und x_0}-\rho_0(r)dr\|_{\mathcal I_N} >
 c'N^{\beta+\alpha_1-1}\Big] \le cc'  N^{1-2\alpha_1}
  \end{equation}

For the last term in \eqref{DDF.56} we bound
 	 \begin{equation*}
\|\pi_{\und x^{\delta,-}_0}-N\rho^{\delta,-}(dr)\|_{\mathcal I_N} \le
\|\pi_{\und x^{\delta,-}_0}-( e^{-\delta}N -N^{\alpha_0})\rho_0(dr)\|_{\mathcal I_N}
+N^{\alpha_0}
   \end{equation*}
because we need to delete from $\pi_{\und x_0}$ the leftmost $ (1-e^{-\delta})N + N^{\alpha_0}$ particles, while the mass to delete from $N\rho_0(dr)$ is
$ (1-e^{-\delta})N$.
By Proposition \ref{lemmaDDF7.1} and \eqref{DDF.19.0}
  \begin{equation*}
  P^{(N)}\Big[
  \|
  \pi_{\und x^{\delta,-}_0}-( e^{-\delta}N -N^{\alpha_0})\rho_0(dr)
  \|_{\mathcal I_N} >
 c'N^{\beta+\alpha_1}\Big] \le cc'  N^{1-2\alpha_1}
  \end{equation*}
therefore
 \begin{equation*}
  P^{(N)}\Big[
  \|\pi_{\und x^{\delta,-}_0}-N\rho^{\delta,-}(dr)\|_{\mathcal I_N} >
 c'N^{\beta+\alpha_1}+N^{\alpha_0}\Big] \le cc'  N^{1-2\alpha_1}
  \end{equation*}
With the choice $\alpha_0=\alpha_1$, $c'N^{\beta+\alpha_1}+N^{\alpha_0}
\le c''' N^{\beta+\alpha_1}$ hence
\begin{equation*}
  P^{(N)}\Big[
  \|\frac 1N \pi_{\und x^{\delta,-}_0}-\rho^{\delta,-}(dr)\|_{\mathcal I_N} >
 c''' N^{\beta+\alpha_1-1}\Big] \le cc'  N^{1-2\alpha_1}
  \end{equation*}

\smallskip

It remains to consider the second term on the right hand side of \eqref{DDF.56}.
We fix a time interval $(k\delta,(k+1)\delta)$,
call $\nu_t(dr)= N \rho^{\delta,\pm}(r,k\delta+t) dr$, $\und y_0 = \und x^{\delta,\pm}_{k}$
and
	 \begin{equation}
	 \label{DD5.19}
\la^{\und y_0}_t(dr):= \int P^{\und y_0}_t(d\und y)\pi_{\und y}(dr)
   \end{equation}
By the triangular inequality we have:
	 \begin{equation}
	 \label{DD5.23}
\| \pi^{\und y_0}_{\und y(t)}-\nu_t\|_{\mathcal I'_N} \le \| \pi^{\und y_0}_{\und y(t)}-\la^{\und y_0}_t\|_{\mathcal I'_N}+ \| \la^{\und y_0}_t-\nu_t\|_{\mathcal I'_N}
   \end{equation}

   \smallskip
   \noindent
   {\bf Bound of the first term in \eqref{DD5.23}}.

 By \eqref{DD5.24}
	 \begin{equation*}
	 \label{DD5.28}
P_t^{\und y_0}\Big[ \bigcap_{I\in \mathcal I'_N} \{| \pi^{\und y_0}_{\und y(t)}(I)-\la^{\und y_0}_t(I)| < N^{\alpha_1} | \}\ge 1- c' N^{\beta+b}cN^{1-2\alpha_1}
   \end{equation*}
because the cardinality $|\mathcal I'_N|$ of $\mathcal I'_N$ is bounded by  $c'N^{\beta+b}$.
Then:
	 \begin{equation*}
	 \label{DD5.29}
P_t^{\und y(0)}\Big[ \| \pi^{\und y_0}_{\und y(t)}-\la^{\und y_0}_t\|_{\mathcal I'_N}\le
c' N^{\beta+b}  N^{\alpha_1}\Big] \ge 1- c'  N^{\beta+b}cN^{1-2\alpha_1}
   \end{equation*}
 By our choice \eqref{DD5.30} of the parameters the last term vanishes as $N\to \infty$.

   \smallskip
   \noindent
   {\bf Bound of the second term in \eqref{DD5.23}}.   By the triangular inequality
 	 \begin{equation}
	 \label{DD5.32}
 \| \la^{\und y_0}_t-\nu_t\|_{\mathcal I'_N} \le
 \| \la^{\und y_0}_t-\la^{\und y'_0}_t\|_{\mathcal I'_N}+\| \nu'_t-\nu_t\|_{\mathcal I'_N} +\| \la^{\und y'_0}_t-\nu'_t\|_{\mathcal I'_N}
   \end{equation}
where   $\und y'_0$ is obtained from $\und y_0$ by shifting each
$x \in \und y_0$ to the center $x_I$ of the set $I\in \mathcal I_N$ where $x$ belongs. Analogously
   $$
   \nu'_t(dr) = \sum_{I\in \mathcal I_N} \nu_0(I) e^{L t}(x_I,r)dr
   = \sum_{I\in \mathcal I_N} \nu_0(I) \la^{x_I }_t(dr)
   $$
   By \eqref{DD5.22} and Theorem \ref{thmDD5.1}
   	 \begin{equation*}
 \| \la^{\und y_0}_t-\la^{\und y'_0}_t\|_{\mathcal I'_N}\le
 \sum_{x\in \und y_0}  \| \la^{x}_t-\la^{x_I}_t\|_{\mathcal I'_N}
 \le \sum_{x\in \und y_0} \int_{\mathbb R} dy\, |e^{Lt}(x,y)-
 e^{Lt}(x_I,y)|
   \end{equation*}
where given $x \in \und y_0$ $x_I$ is the center of the set $I\in \mathcal I_N$ where $x$ belongs.  By \eqref{DD3.222}
     	 \begin{equation*}
 \| \la^{\und y_0}_t-\la^{\und y'_0}_t\|_{\mathcal I'_N}
 \le   N c \delta^{-1/2} N^{-\beta}
   \end{equation*}
  because  $ |x-x_I| \le N^{-\beta}$.  In an analogous way we prove that
        	 \begin{equation*}
	 \label{DD5.35}
 \| \nu'_t-\nu_t\|_{\mathcal I'_N}
 \le  N c \delta^{-1/2} N^{-\beta}
   \end{equation*}
 The last term in \eqref{DD5.32} is bounded by
 	 \begin{equation*}
	 \label{DD5.36}
 \| \la^{\und y'_0}_t-\nu'_t\|_{\mathcal I'_N} \le \sum_{J\in \mathcal I_N}
  \| m_J\la^{x_J}_t-\nu_0(J)\la^{x_J}_t\|_{\mathcal I'_N}
  = \sum_{J\in \mathcal I_N} | m_J -\nu_0(J)|
  \|\la^{x_J}_t\|_{\mathcal I'_N}
   \end{equation*}
 $m_J$ the number of particles in $\und y_0$ which are in $J$.  Hence
   	 \begin{equation*}
	 \label{DD5.37}
 \| \la^{\und y'_0}_t-\nu'_t\|_{\mathcal I'_N} \le \sum_{J\in \mathcal I_N}
 e^t | m_J -\nu_0(J) | =e^t \| \pi_{\und y_0}-\nu_0\|_{\mathcal I_N}
   \end{equation*}
In conclusion
	 \begin{equation}
	 \label{DD5.38}
 \| \la^{\und y_0}_t-\nu_t\|_{\mathcal I'_N} \le 2 c  \delta^{-1/2}   N^{1-\beta}
 + e^t \| \pi_{\und y_0}-\nu_0\|_{\mathcal I_N}
   \end{equation}

Thus by \eqref{DD5.23}, \eqref{DD5.29} and \eqref{DD5.38}:
	 \begin{eqnarray}
	 \label{DDF.5tris}
&& \sum_{k=1}^{2^n}\{P^{(N)}[
(\mathcal X^{+}_2(k))^c] +P^{(N)}[
(\mathcal Z^{-}_2(k))^c]\} \le 2^n 2c' cN^{1+\beta+b-2\alpha_1}
   \end{eqnarray}

\appendix

\section{Probability estimates}

The number of particles $n(t)=|\und y(t)|$ in the basic process is a Markov process with state space $\mathbb N$ and generator
	 \begin{equation}
	 \label{DD55.1}
  K f(n)  = n \Big( f(n+1) - f(n) \Big)
   \end{equation}
   Denote by $P^{n_0}$ the law of $n(t)$ starting from $n_0$ and by $E^{n_0}$ its expectation.

\begin{lemma}
\label{lemmaDD55.1}

For any positive integer $k$ there is $c_k$ so that
	 \begin{equation}
	 \label{DD55.2}
 E^{n_0}[n(t)] \le  e^{c_k t}  n_0
   \end{equation}

\end{lemma}

\medskip

\noindent
{\bf Proof.}  Call $h_k(t)$ the left hand side of \eqref{DD55.2}, then
	 \begin{equation*}
 \frac{d}{dt}h_k(t) =  E^{n_0}[n(t)\{ (n(t)+1)^k - n(t)^k]\le k
  E^{n_0}[n(t) k (n(t)+1)^{k-1}] \le
 c_k h_k(t)
   \end{equation*}
   $c_k= k 2^k$.
\qed

\begin{lemma}
\label{lemmaDD5.4}
Suppose that $\und y_0$ is a configuration with $\le N$ particles, let $\alpha_1 \in (\frac 12,1)$, $T>0$. Call $\la_t^{\und y_0}(dr)= \sum_{x \in \und y_0}\la_t^{x}(dr)$, see \eqref{DD5.1}.
Then there is $c$ so that for any   $t \le \delta$ and any interval $I\in \mathcal I_N$:
	 \begin{equation}
	 \label{DD5.24}
P_t^{\und y_0}\Big[ | \pi^{\und y_0}_{\und y(t)}(I)-\la^{\und y_0}_t(I)| \ge N^{\alpha_1} | \le c N^{1-2\alpha_1}
   \end{equation}
	 \begin{equation}
	 \label{DD5.25}
P_t^{\und y_0}\Big[ | \pi^{\und y_0}_{\und y(t)}(\mathbb R)-\la^{\und y_0}_t(\mathbb R)| \ge N^{\alpha_1} \Big] \le c \delta N^{1-2\alpha_1}
   \end{equation}

\end{lemma}

\noindent
{\bf Proof.}
By \eqref{DD5.22}
	 \begin{equation}
	 \label{DD5.26}
 |\pi^{\und y_0}_{\und y(t)}(I)-\la^{\und y_0}_t(I)| = |
\sum_{x\in \und y_0} \{\pi^{x}_{\und y(t)}(I)-\la^{x}_t(I)\}|
   \end{equation}
The variables $\{\pi^{x}_{\und y(t)}(I)-\la^{x}_t(I)\}$ are mutually independent and centered then by the Chebishev inequality with the squares we get:
	 \begin{equation*}
P_t^{\und y_0}\Big[  |\pi^{\und y_0}_{\und y(t)}(I)-\la^{\und y_0}_t(I)|\ge N^{\alpha_1} | \Big]\le N^{-2\alpha_1}\sum_{x \in \und y(0)} \{ \int P^x_t(d\und y) \pi_{\und y}(I) ^2 - \la^x_t[I]^2\}
   \end{equation*}
and \eqref{DD5.24} follows from \eqref{DD55.2}.  The left hand side of \eqref{DD5.25} is bounded by
	 \begin{equation*}
\le N^{-2\alpha_1}\sum_{x \in \und y(0)} \{ \int P^x_t(d\und y) \pi_{\und y}(\mathbb R) ^2 - \la^x_t[\mathbb R]^2\}
   \end{equation*}
We have
  $$
   \int P^x_t(d\und y) \pi_{\und y}(\mathbb R) ^2= 2e^{2t}-e^{t}, \quad \la^x_t[\mathbb R]^2= e^{2t}
   $$
hence  \eqref{DD5.25} recalling that the cardinality of $ \und y(0)$ is $\le N$.
  \qed

\begin{theorem}
\label{thmDD5.3}
Let $\und y_0$ be a configuration with $\le N$ particles all in the support $[-A,A]$ of $\rho_0$.
Fix $T>0$ and $b>0$ then there is $c_{b,T}$ so that
	 \begin{equation}
    \label{DD5.13}
P^{\und y_0}\Big[ \und y(t) \cap [-N^b, N^b]^c = \emptyset, 0\le t \le T\Big] \ge 1 - c_{b,T} N^{-1}
 \end{equation}

\end{theorem}

\noindent
{\bf Proof.} We will only prove that
	 \begin{equation}
    \label{DD5.14}
P^{\und y_0}\Big[ \und y(t) \cap [N^b,\infty) = \emptyset, 0\le t \le T\Big] \ge 1 - c N^{-1}
 \end{equation}
as the bound for the probability of $ \und y(t) \cap (-\infty,-N^b] = \emptyset$
 is similar.
 We have
	 \begin{eqnarray}
    \label{DD5.15}
&&P^{\und y_0}\Big[\text{there is $t\in [0,T]$}: \;\und y(t) \cap [N^b,\infty) \ne \emptyset \Big]\nn\\&&\hskip2cm \le  N   P^{\{A\}}\Big[\text{there is $t\in [0,T]$}: \;\und y(t) \cap [N^b,\infty) \ne \emptyset \Big]
 \end{eqnarray}
where $P^{\{A\}}$ is the law of  $\und y(t)$ starting from a single particle at position $A$, recall that
$[-A,A]$ is the support of the initial density $\rho_0$ and that
$[-\xi,\xi]$ is the support of  $p(0,x)$. Let $\und z(t)$ be the process of branching brownian particles where at any branching time a new particle is put at $x+\xi$ (if generated by a particle at $x$) while all the particles except the new one are shifted to the right by $\xi$.  Thus $\und z(t)= \und x^0(t)+ n(t)\xi$ where $n(t)$ is the number of branching times till $t$ and $\und x^0(t)$ is the process where the branching is local ($p(x,y)= \delta(y-x)$) and no particle  is killed.   Calling $\mathbb P^{0}$ the law of the process $\und x^0(t)$ starting from a single particle at position $0$,
	 \begin{eqnarray}
    \label{DD5.16}
&&NP^{\{A\}}\Big[\text{there is $t\in [0,T]$}: \;\und y(t) \cap [N^b,\infty) \ne \emptyset \Big]\nn\\&&\hskip2cm \le  N \mathbb P^{0}\Big[\text{there is $t\in [0,T]$}: \;\und x^0(t) \cap [N^b-n(t)\xi-A,\infty) \ne \emptyset \Big]
 \end{eqnarray}
By Lemma \ref{lemmaDD55.1} for any integer $k$ there is $c_k$ so that,
	 \begin{eqnarray}
    \label{DD5.17}
&&\mathbb P^{0}\Big[ n(T) \ge N^a \Big] \le c_k N^{-ak}
 \end{eqnarray}
which will be used with  $a=b/2$.  \eqref{DD5.17} follows from the fact that for any $k$ the expectation
 $\mathbb E^{0} \Big[ n(T) ^k \Big]$ is bounded. We choose $k$ so that $ak >2$, then the right hand side of \eqref{DD5.16} is bounded by
 	 \begin{equation}
    \label{DD5.18}
\le  c_k N^{-1} +
N \mathbb P^{0}\Big[n(t) <N^a;\;\text{there is $t\in [0,T]$}: \;\und x^0(t) \cap [N^b-n(t)\xi-A,\infty) \ne \emptyset \Big]
 \end{equation}
 Let $n(t) = k$ and $\mathcal T$ one of the trees obtained from the branching history
 with $k$ outputs.  To each branch of the tree we associate a particle, the law of its motion is that of a brownian motion $B(t)$ starting from 0.  Since there are at most $N^a$ branches, the second term in \eqref{DD5.18} is bounded by $\le
N^{1+a} P \big[ \max_{0\le t \le T},|B(t)| \ge N^b-N^a\xi-A\big]$
 hence \eqref{DD5.13}. \qed

\begin{corollary}
\label{thmDD5.3.1}

With the notation of Theorem \ref{thmDD5.3}
	 \begin{equation}
    \label{DD5.13.1}
P\Big[ \und x^{\delta,\pm}(t) \cap [-N^b, N^b]^c = \emptyset, 0\le t \le T\Big] \ge 1 - c_{b,T} N^{-1}
 \end{equation}

\end{corollary}

\noindent
{\bf Proof.}  The processes $ \und x^{\delta,\pm}(t) $ can be realized as subsets of the process $\und y(t)$, hence \eqref{DD5.13.1}. \qed

\def\at{,\ }
\def\email#1{{\tt #1}}
{
Anna De Masi\at
              Universit\`a di L'Aquila, 67100 L'Aquila, Italy 
                            \email{demasi@univaq.it}           

  Pablo A. Ferrari\at
              Universidad de Buenos Aires, DM-FCEN, 1428 Buenos Aires, Argentina
              \email{pferrari@dm.uba.ar}

          Errico Presutti\at
              Gran Sasso Science Institute, 67100 L'Aquila, Italy\\
              \email{errico.presutti@gmail.com} 

           Nahuel Soprano-Loto\at
              Universidad de Buenos Aires, DM-FCEN, 1428 Buenos Aires, Argentina
              \email{sopranoloto@gmail.com}
}

\end{document}